\documentclass[preprint]
{iucr}        

\journalcode{A}              

\usepackage[usenames]{color}
\usepackage{graphicx}
\usepackage[english]{babel}
\usepackage[dvipsnames]{xcolor}
\usepackage{amsthm}
\usepackage{amsmath}
\newcommand{\seqnum}[1]
{{\underline{#1}}}

\newcommand{\url}[1]{\underline{#1}}
\newcommand{\beql}[1]{\begin{equation}\label{#1}}
\newcommand{\eeq}{\end{equation}}
\newcommand{\eqn}[1]{(\ref{#1})}

\newtheorem{thm}{Theorem}{\bfseries}{\itshape}
{\bfseries}{\itshape}
{\bfseries}{\itshape}
{\bfseries}{\itshape}
{\bfseries}{\itshape}

\newcommand{\ldot}{\!.}

\begin{document}           



\title{A Coloring Book Approach to Finding Coordination Sequences}
\shorttitle{Finding Coordination Sequences}


\author[a]{C.}{Goodman-Strauss}
\cauthor[b]{N.J.A.}{Sloane}{njasloane@gmail.com} 

\aff[a]{SCEN 303, Univ. Arkansas, \city{Fayetteville}, AR 72701, \country{USA}}
\aff[b]{The OEIS Foundation Inc., 11 So. Adelaide Ave., \city{Highland Park} NJ 08904 \country{USA}}
\shortauthor{Goodman-Strauss and Sloane}

[Publication information:
 Chaim Goodman-Strauss and N. J. A. Sloane,
  A Coloring Book Approach to Finding Coordination Sequences,
  Acta Cryst. Sect.A: Foundations and Advances, 2019, A75:121-134;
https://doi.org/10.1107/S2053273318014481]

\maketitle                        

\begin{synopsis}
This article presents a simple method for finding formulas for coordination sequences, based on coloring the underlying graph according to certain rules. It is illustrated by
applying it to several uniform tilings and their duals. 
\end{synopsis}

\begin{abstract}
An elementary method is described for finding the coordination sequences
for a tiling, based on coloring the underlying graph. 
The first application is to the two kinds of vertices (tetravalent and trivalent) in the
Cairo (or dual-$3^2 \ldot 4.3.4$) tiling.
The coordination sequence for a tetravalent
vertex turns out, surprisingly, to be $1, 4, 8 ,12, 16, \ldots$,
the same as for a vertex in the familiar square (or $4^4$) tiling.
The authors thought that such a simple fact should have a simple proof,
and this article is the result.
The method is also used to obtain coordination sequences for the $3^2 \ldot 4.3.4$,
$3.4.6.4$, $4.8^2$,  $3.12^2$, and $3^4 \ldot 6$ uniform tilings,
and the snub-$632$ tiling.
In several cases the results
provide proofs for previously conjectured formulas.
\end{abstract}



\section{Introduction}\label{Sec1}

The ``Cairo tiling'' (Figure~\ref{FigCairo}) has many names,
as we will see in Section~\ref{SecCairo}.  In particular, it is the dual of the  snub version of the familiar square tiling.  It has two kinds of vertices (i.e. orbits of vertices under its symmetry group) --- tetravalent and trivalent.

It therefore has two coordination sequences (CS's), giving the numbers  $a(n)$ of vertices that are at each distance $n$ from a base vertex $P$, as measured in the graph of edges and vertices  in the tiling: those based at either a 
trivalent or tetravalent vertex. These coordination sequences may be read off from the contour lines, of equal distance from the base vertex $P$ (Figure~\ref{FigTet}, left and right, respectively). However, although  contour lines are structured and can be used to compute CS's,  they  are  unwieldy to construct and analyze.
It is usually easy enough to find the first few terms of the CS, but our goal is to find a formula, recurrence, or generating function for the sequence.

This article was motivated by our recent discovery that the
coordination sequence  for a tetravalent vertex
in the Cairo tiling appeared to be the same as that for the square grid. We thought that such a simple fact should have a simple proof,
and this article is the result.

The {\em coloring book approach}, 
described in \S\ref{SecM}, is an elementary means of calculating coordination sequences, based on coloring the underlying graph with ``trunks and branches'' and finding a recurrence for the number
of vertices at a given distance from the base point.
We have to verify that a desired local structure propagates, and that our colored trees  
do give the correct distance to the base vertex.
In \S\ref{SecGrid} we  illustrate the method by applying it to an elementary case, 
the $4^4$ square grid tiling. 

The method works quite well in many cases and is at least helpful in others:
in Sections \S\ref{SecTet} and \S\ref{SecTri} we
deal with the tetravalent  and trivalent vertices  in the Cairo
tiling and in \S\ref{SecDual} with its dual, the uniform (or Archimedean) $3^2 \ldot 4.3.4$ tiling. 
We then apply the method to obtain coordination sequences for four other uniform tilings,
$3.4.6.4$ (\S\ref{Sec3464}),
$4.8^2$ (\S\ref{Sec488}),  $3.12^2$ (\S\ref{Sec31212}),
$3^4 \ldot 6$ (\S\ref{Sec33336}),
and to the dual of the $3^4 \ldot 6$ (\S\ref{SecSnub}) tiling.

Starting in \S\ref{Sec31212} we  must rely on a more subtle analysis, but find that the coloring book method at least renders our calculations somewhat more transparent.

There are of course many works that discuss more sophisticated methods
for calculating coordination sequences, in both the 
crystallographic and mathematical literature, such as
Baake \& Grimm (1997),
Bacher \& de la Harpe (2018), 
Bacher, de la Harpe \& Venkov, B. (1997),
Conway \& Sloane (1997),
de la Harpe (2000),
Eon (2002, 2004, 2007, 2013, 2016, 2018), 
Grosse-Kunstleve, Brunner \& Sloane (1996),
O'Keeffe (1995), and
O'Keeffe and Hyde (1980).
For a uniform tiling, where the symmetry group acts transitively
on the vertices, there is an alternative method for finding CS's based on 
Cayley diagrams and the Knuth-Bendix algorithm, and using
the computer algebra system Magma 
(Bosma, Cannon \& Playoust, 1997).
This approach is briefly described in \S\ref{SecCD}.
For uniform tilings it is known (Benson, 1983) that that the CS has a rational generating function.  This is also known to be true for other classes of tilings (see the references mentioned at the start of this paragraph). But as far as we know the general conjecture that the CS of a periodic tiling of $n$-dimensional Euclidean space by polytopes always has a rational generating function is still open, even in two dimensions. In this regard, the work of Zhuravlev (2002) and others on the limiting shape of the contour lines in a two-dimensional tiling may be relevant.

A more traditional way to calculate the CS by hand is to draw `contour lines'  or 'level curves' that
connect the points at the same distance from $P$.
These lines are usually overlaid on top of the tiling.
The resulting picture can get very complicated (see Fig.~\ref{FigTet}), and this approach is
usually only successful for simple tilings or for finding just
the first few terms of the CS.

``Regular production systems'' (Goodman-Strauss, 2009)
can be used to give
a formal model of growth along a front, 
enabling  generalized contour lines to be described as languages, and underlies some of our thinking here.

The computer program {\em ToposPro}
(Blatov,  Shevchenko \& Proserpio, 2014)
makes it easy to compute the initial terms of the coordination sequences 
for a large number of tilings, nets (both two- and three-dimensional),
crystal structures, etc. 

Besides the obvious application of coordination sequences for
estimating the density of points in a tiling, another use is for identifying
which tiling or net is being studied. This is especially useful 
when dealing with three-dimensional structures, as in 
(Grosse-Kunstleve, Brunner \& Sloane,1996).
Another application of our coloring-book approach is for finding 
labels for the vertices in a graph, as we mention at the end of \S\ref{SecM}.

For any undefined terms about tilings, see the classic work 
Gr\"{u}nbaum  \& Shephard (1987)
or the article 
Gr\"{u}nbaum  \& Shephard (1977).


\section{The Cairo tiling}\label{SecCairo}

\begin{figure}
\caption{A portion of the Cairo tiling.} 
\centerline{\includegraphics[angle=0, width=.9\textwidth]{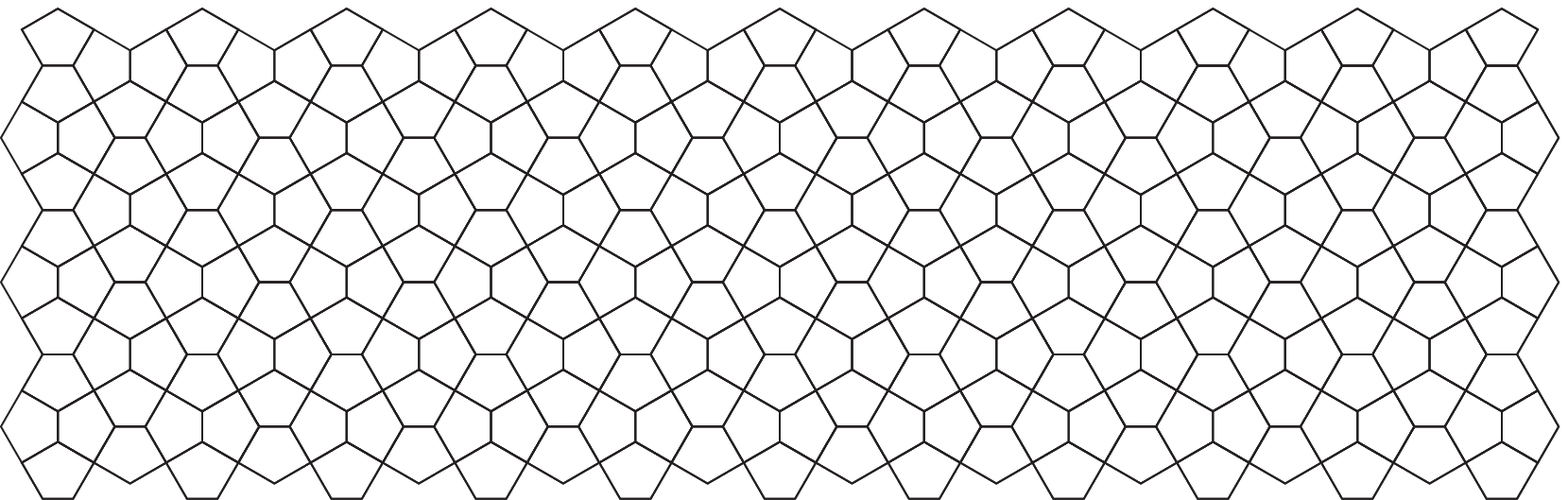}}
\label{FigCairo}
\end{figure}

\begin{figure} 
\centerline{
\includegraphics[width=.4\textwidth]{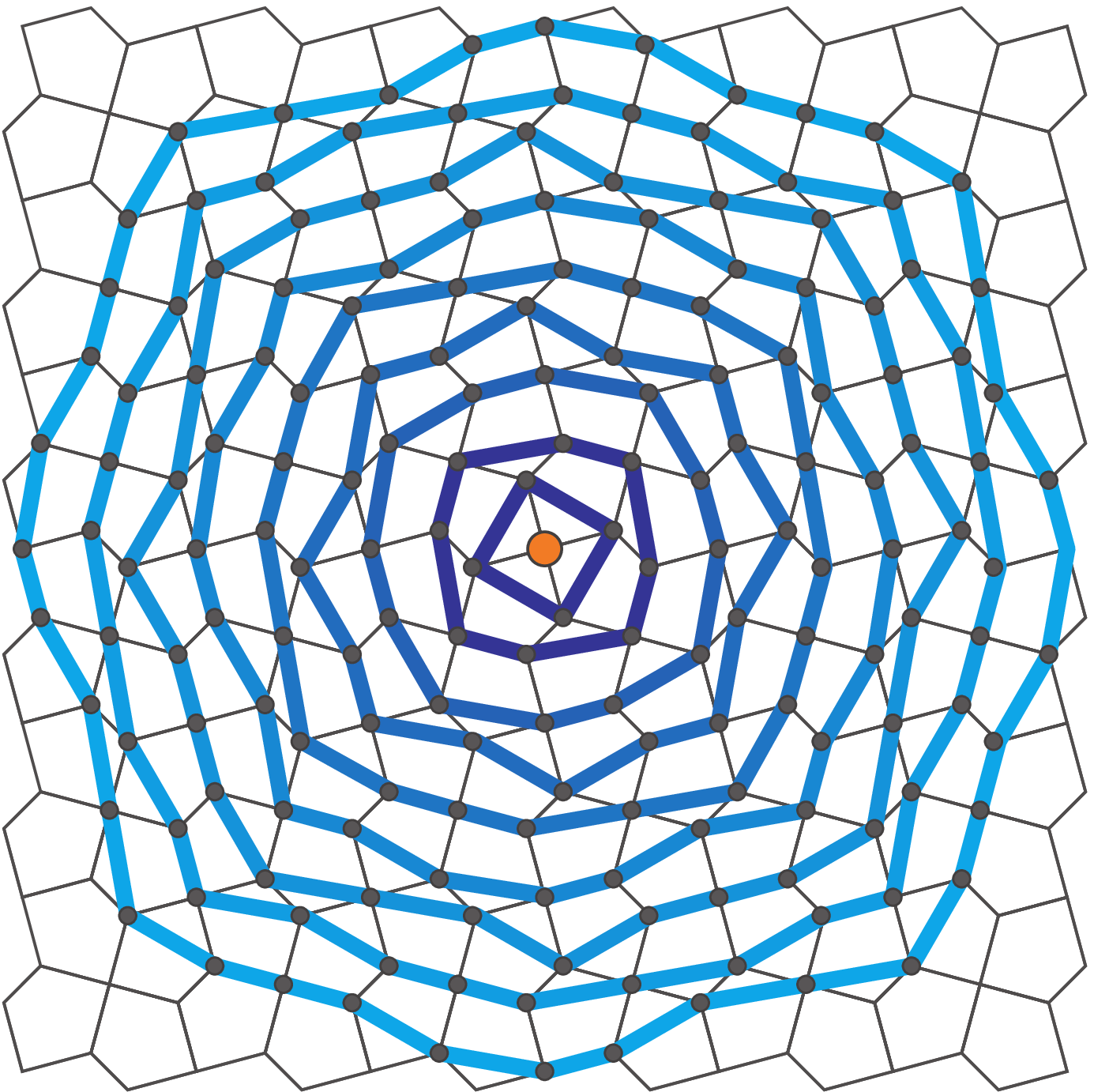}\hspace{.1\textwidth}
\includegraphics[width=.4\textwidth]{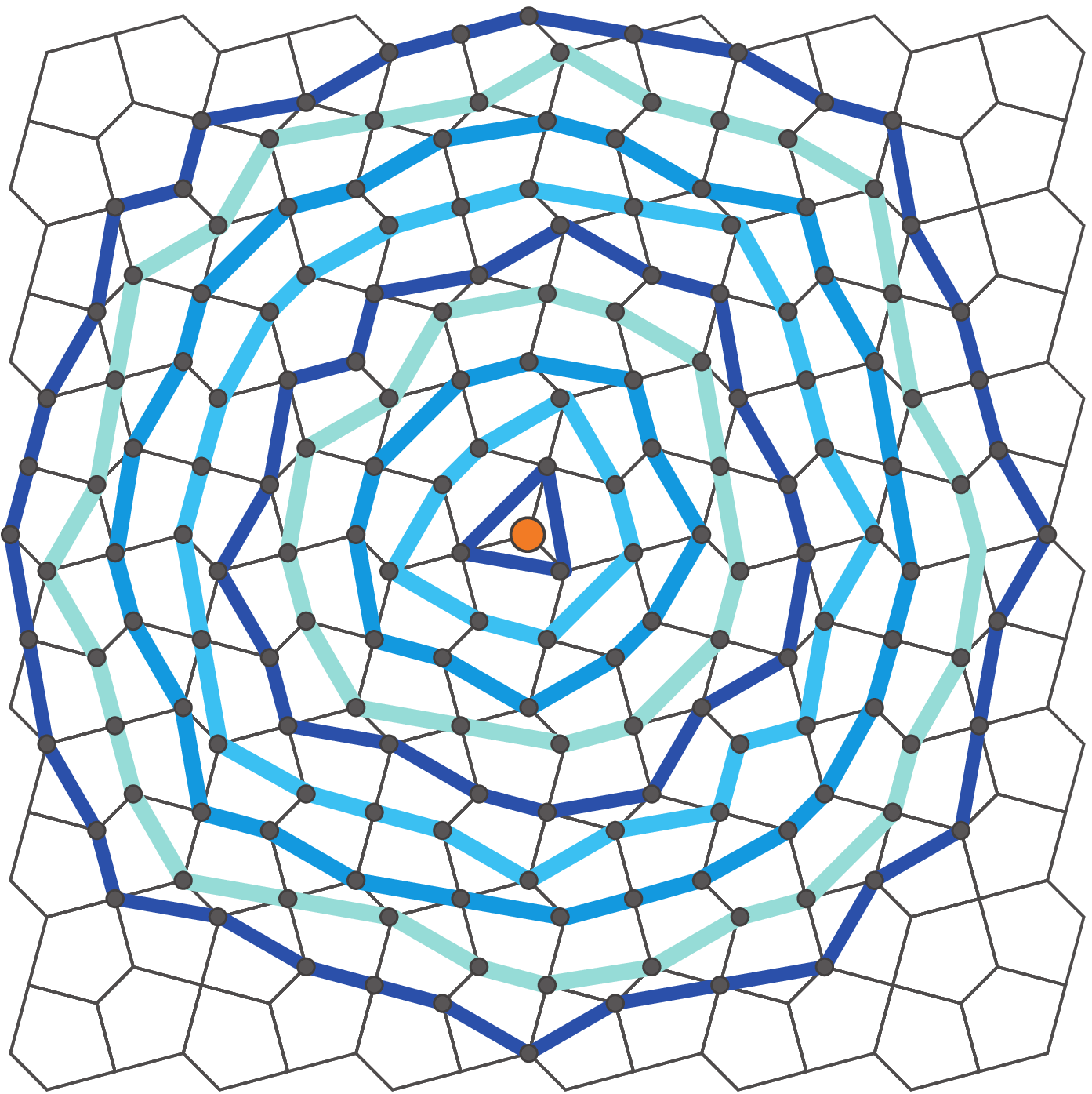}}
\caption{ The Cairo tiling has two kinds of vertices (i.e., orbits of vertices under the $4*2$ symmetry group of the tiling) --- tetravalent and  trivalent. At left, the contour lines of constant distance from a central tetravalent vertex, from which we can read off the first several terms of its CS. At right, contours centered at a trivalent vertex (\S\ref{SecTri}). 
These contour lines have a recursive and analyzable structure,
but  we find that the {\em 
coloring book} approach, shown in Figures~\ref{FigTet2} and~\ref{FigTri2}, is much simpler, allowing more easily verifiable calculations. 
}
\label{FigTet}\label{FigTri}
\end{figure}

The Cairo tiling is shown in Fig. \ref{FigCairo}.
This beautiful tiling has many names. It has also been called
the Cairo pentagonal tiling, 
the MacMahon net (O'Keeffe  \& Hyde, 1980),
the mcm net (O'Keeffe et al., 2008),
the dual of the $3^2 \dot 4.3.4$ tiling 
(Gr\"{u}nbaum  \& Shephard, 1987, pp.~63, 96, 480 (Fig.~$P_5$-24))
the dual snub quadrille tiling,  
or the dual snub square tiling  
(Conway,  Burgiel \& Goodman-Strauss, 2008, pp.~263, 288).
\bigskip
\hrule
\bigskip
Figure~\ref{FigCairo} around here
\bigskip
\hrule
\bigskip

We will refer to it simply as the Cairo tiling.
There is only one shape of tile, an irregular pentagon, which may be 
varied somewhat.\footnote{The tiling may be modified   
without affecting its combinatorics and  coordination sequences, so long as we the topology of the orbifold graph is preserved (cf.~Conway et al., 2008).}

The tiling is named from its use in Cairo, 
where this pentagonal tile has been mass-produced 
since at least the 1950's and is prominent around the city. 

\bigskip
\hrule
\bigskip
Figure~\ref{FigTri} around here
\bigskip
\hrule
\bigskip

The CS with respect to a tetravalent vertex in the Cairo tiling begins 
\beql{EqGrid}
1, 4, 8, 12, 16, 20, 24, 28, 32, \ldots,
\eeq
which suggests that it is the same as the CS of
the familiar $4^4$ square grid 
(sequence \seqnum{A008574}\footnote{Six-digit numbers prefixed by A refer to entries in the On-Line Encyclopedia of Integer Sequences, or OEIS.} 
in the OEIS) . We will
show in \S\ref{SecTet} that this is true, by proving:
\begin{thm}\label{ThTet}
The coordination sequence with respect to a tetravalent vertex in
the Cairo tiling is given by $a(0)=1$, $a(n)=4n$ for $n \ge 1$.
\end{thm}

The CS with respect to a trivalent vertex begins 
\beql{EqA296368}
1, 3, 8, 12, 15, 20, 25, 28, 31, 36, 41, 44, 47, 52, 57, 60, 63, 68, 73, 76, \ldots,
\eeq
which has now been added to the OEIS as sequence \seqnum{A296368}. In \S\ref{SecTri} we will prove:
\begin{thm}\label{ThTri}
The coordination sequence with respect to a trivalent vertex in
the Cairo tiling is given by $a(0)=1, a(1)=3, a(2)=8$, and,
for $n \ge 3$,
\begin{align}\label{EqTrivalent}
  a(n) & ~=~ 4n ~~~~~  \mbox{~if~$n$~is~odd} \,,  \nonumber \\
         & ~=~ 4n-1   \mbox{~if~} n \equiv 0 \pmod{4} \,,  \nonumber \\
         & ~=~ 4n+1   \mbox{~if~} n \equiv 2 \pmod{4} \,.          
\end{align}
\end{thm}

\begin{figure} 
\centerline{\includegraphics[width=.4\textwidth]{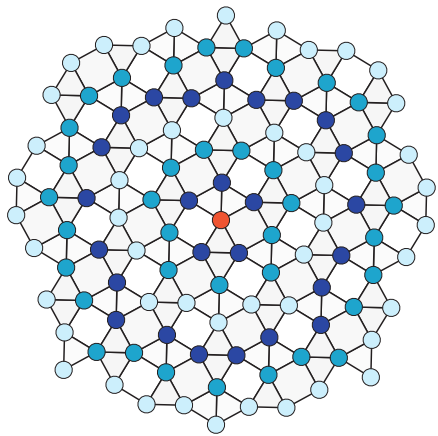}}
\caption{A portion of the $3^2 \ldot 4.3.4$ uniform tiling, the dual to the Cairo structure.
Distances to $P$ (which is the red vertex in the center of the figure) are shown in various colors,
illustrating terms $a(0)$, \ldots, $a(6)$ of the coordination sequence,  \seqnum{A219529}. In Figure~\ref{FigDualb}, the count is clarified by  a trunks and branches structure.}
\label{FigDual}
\end{figure}

\bigskip
\hrule
\bigskip
Figure~\ref{FigDual} around here
\bigskip
\hrule
\bigskip

The dual of the Cairo tiling is the uniform (or Archimedean) tiling $3^2 \ldot 4.3.4$ shown in Fig.~\ref{FigDual}. 
Now there is only one kind of vertex, with valency $5$, and  the coordination sequence begins
\beql{EqDual}
1, 5, 11, 16, 21, 27, 32, 37, 43, 48, 53, 59, 64, 69, 75, 80, 85, 91, 96, 101, \ldots,
\eeq
given in \seqnum{A219529}.  That entry has a long-standing conjecture 
that $a(n) = \lfloor  \frac{16n+1}{3} \rfloor$, and we will establish this in \S\ref{SecDual}
by proving:
\begin{thm}\label{ThDual}
The coordination sequence with respect to a vertex in
the $3^2 \ldot 4.3.4$ tiling is given by $a(0)=1$,  
\begin{align}\label{EqDuala}
  a(3k) & ~=~ 16k, ~   k \ge 1 \,,  \nonumber \\
  a(3k+1)  & ~=~ 16k+5 , ~ k \ge 0 \,,  \nonumber \\
 a(3k+2)   & ~=~ 16k+11, ~ k \ge 0  \,.          
  \end{align}
\end{thm}


\section{The coloring book approach to finding coordination sequences}\label{SecM}

We start with a more precise statement of the problem.
Let $T$ be a periodic tiling of the plane by polygonal tiles.
The graph $G = G(T)$ of the tiling has a vertex for each point of the plane
where three or more tiles meet, and an edge between two vertices
if two tiles share a boundary along the line joining the
corresponding points.

We assume the tiling is such that $G$ is a connected graph:
$G$ is thus a connected, periodic, planar graph with all vertices
of valency at least $3$, and (since the tiles are polygons) with no parallel edges.
The coloring book method could be applied
to any  graph of this type, not just one arising from a tiling.

The distance $d(Q,R)$ between vertices $Q$, $R$ in $G$ is defined to be  the number of
edges in the shortest path joining them. The coordination sequence of $G$ with
respect to a vertex $P \in G$ is the sequence
$\{a(n): n = 0, 1, 2, \ldots\}$
where $a(n)$ is the number of vertices $Q \in G$ with $d(Q,P)=n$.
We refer to $P$ as the {\em base vertex}.
For a periodic tiling there are only a finite number
of different choices for the base vertex, and our goal
is to find the coordination sequence with respect to a base vertex
of each possible type.

Our method for finding the coordination
sequence with respect to a base vertex $P$ is to try to find
a subgraph $H$ of $G$
with the following properties.
\begin{list}{}{\setlength{\itemsep}{0.02in}}
\item[(i)] $H$ is a connected graph that passes through every vertex of $G$, and
\item[(ii)] for any vertex $Q$, every path in $H$ from $Q$ to $P$ has the minimal
possible length, $d(Q,P)$.
\end{list}

We also want $H$ to have three further, less well-defined, properties.

\begin{list}{}{\setlength{\itemsep}{0.02in}}
\item[(iii)] $H$ should be essentially a tree (in
the sense of graph theory), and more precisely should consist of a finite
number of ``trunks'' (in the arboreal sense), which are disjoint paths that originate at $P$,
together with infinitely many ``branches'', which are also disjoint paths
and originate at trunk vertices. However, on occasion we will allow
the trunks to have ``burls'' (i.e., bulges, or loops inside the trunk), and both trunks and branches may have ``twigs'' (typically small subtrees with just a few edges) 
growing from them.

In our figures, trunks will usually be colored
\textcolor{blue}{blue}) and branches  \textcolor{OliveGreen}{green}.
A glance at some of the figures below will illustrate our arboreal
terminology.  In Figs.~\ref{FigGrid} and \ref{FigTet2} (left)
$H$ actually {\em is} a tree, with simple trunks and branches.  
In Fig.~\ref{FigTri2} (left) $H$ is still
a tree, but there are twigs (single edges) between the two
parallel trunks.
In Figs.~\ref{FigDualb} (left and right) two
of the trunks have burls (loops of lengths $4$ and $6$, respectively),
and so are not pure trees in the mathematical sense.
In Fig.~\ref{Fig33336Coloring} both the trunks
and branches have twigs attached.

\item[(iv)] There should an easy way to check that Property (ii) holds, i.e., that there
are no shortcuts to $P$ that take a path that is not part of $H$.
\item [(v)] And towards this end, $H$ should have some sort of regular structure that allows us to make inductive arguments.
\end{list}

There is no difficulty in satisfying (i), since any spanning tree rooted at P would do.
But this is not very helpful since there are an infinite number
of distinct spanning trees, and  in any case
requiring $H$ to be a tree in the mathematical sense can make it harder to
satisfy the other conditions.

Assuming that (ii) and (iii) hold, the subgraph $H$ has the following
structure.  Each trunk out of $P$ is an infinite path, and each trunk vertex 
(ignoring the slight complication caused by the
burls) is joined to a unique trunk vertex that is one step further away from $P$.
The branches are infinite paths originating at trunk vertices and
(ignoring the twigs)
each branch vertex is joined to a unique branch vertex that is one step further away 
from $P$. The twigs are finite (and small) subtrees that connect any remaining 
vertices to the closest trunk or branch. If we have carried out (v)  well at all, 
then all these vertices should be easy to count. 

As we have gained experience with the method, we have found the following 
further condition to be helpful for making inductive arguments.
\begin{list}{}{\setlength{\itemsep}{0.02in}}
\item [(vi)] The distance from any vertex of the tiling to
the closest trunk or branch should be bounded. Equivalently,
there should be a constant $c_0$ such that the twigs have length at most $c_0$.
\end{list}
 
An alternative way to describe $H$ is to think
of $G$ as a topographic map, where the heights above sea level of the vertices
are the distances from the base point $P$.
Then $H$ represents a drainage network
that always flows downhill.  In this
model the `trunks' represent major rivers that flow to $P$,
and the `branches' are tributaries that feed into the major rivers.
 
Speaking informally, the subgraph $H$ is usually orthogonal to the contour lines,
just as in polar coordinates, radial lines are orthogonal to the circles.

We have tried a few strategies for finding the subgraph $H$.
In several examples (\S\ref{SecDual}-\S\ref{Sec488}), our human visual system seems to  fill in the proof of Property (ii) instantly. 
But a verifiable means of testing (ii) is essential.
For the Cairo tiling, we can redraw the graph so that (ii) is clear,
as at right in Figures~\ref{FigTet2} and~\ref{FigTri2}. 

As described further in \S\ref{Sec33336}, in later examples we  use  an atlas of patterns, of what locally appear to be trunks and branches, directed across alleged contour lines.
These patterns  propagate outwards from a region about $P$ --- each patch 
extends outwards in a natural way.
Since the alleged contour lines are initially simple  nested  closed curves, 
they continue to be so, and therefore are indeed contour lines. Since the alleged trunks and branches are transverse to the contour lines, they are indeed trunks and branches
satisfying (ii). 

We have found the process of searching for trunks and branches
by drawing with colored pencils on pictures of
tilings to be quite enjoyable.
If readers wish to try this for themselves --- and perhaps
to improve on our constructions --- we encourage them to
download pictures of tilings from the Internet.
There are many excellent web sites.
Galebach's web site (Galebach 2018)
is especially important, as it includes pictures of all $k$-uniform
tilings with $k \le 6$, with over $1000$ tilings.
The Chavey (1989) article and the
Hartley (2018) and Printable Paper (2018) 
web sites have many further pictures, and the
Reticular Chemistry Structure Resource 
or RCSR (O'Keeffe et al. 2008) and 
ToposPro (Blatov et al. 2014) 
databases have thousands more.

In this paper we have been careful to describe our approach as a {\em method}
for finding coordination sequences, rather than an {\em algorithm}.
At present we do not have enough experience with the method
to state it any more formally.  We can point out that periodic two-dimensional  tilings
by polygons fall into two classes: the essential underlying periodic structure
is either rectangular or hexagonal. In the former case one should look
for a trunks and branches structure like that shown in Figures
\ref{FigGrid}, \ref{FigTet2}, \ref{FigTri2},
\ref{Fig3464}, \ref{Fig488},
and in the latter case like that shown in Figure \ref{Fig33336Coloring}, 
\ref{Fig:snubColoring2}. But as one can see by looking at the catalogs of tilings 
mentioned above, there is too great a variety of tilings for us to be any more precise than that.

One's first guess at a trunks
and branches structure  does not always succeed. 
For example, in  the $3.4.6.4$ tiling, shown at right in Fig.~\ref{Fig3464},
a natural first guess is to make the branches roughly  horizontal.
However, this does not provide the shortest paths to the origin --- for that one
needs to use vertical branches that head North-East (in the first quadrant)
and North-West (in the second quadrant), as shown in the figure.

There is sometimes another side-benefit to our approach: it may
provide a way to  assign coordinates to  
the vertices  of the graph. 
If the branches are paths, and
a vertex $Q$ is on a branch that originates at a trunk vertex $R$, then we can label
$Q$ by specifying the trunk, the distance $d(R,P)$,  the branch, and the distance $d(Q,R)$
(with appropriate modifications in case the trunk or branch is not quite a simple path).


\section{The square grid}\label{SecGrid}

\begin{figure} 
\centerline{\includegraphics[width=.45\textwidth]{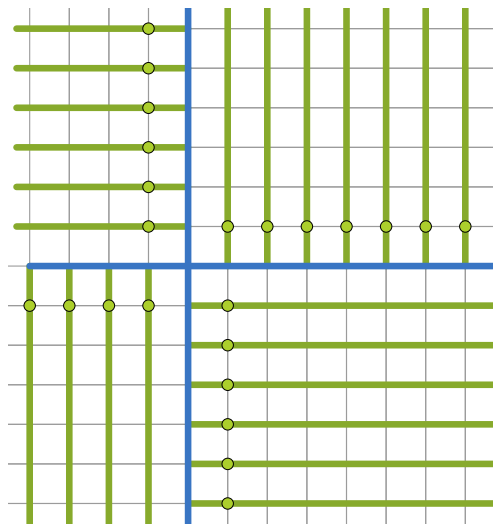}}
\caption{Trunks (blue) and branches (green) for the  familiar square  tiling, $4^4$.
Each vertex is associated (via a blue or green edge) with a unique
vertex that is one step further away from $P$,
except for the vertices on the (blue) trunk,  where branches (green) have sprouted.
As four branches sprout at each distance $n>0$ from $P$,
joining the trunk to the four green vertices,
the  CS satisfies the recurrence  $a(n+1)= a(n)+4$.}
\label{FigGrid}
\end{figure}

Before analyzing the Cairo tiling, we illustrate the coloring book method 
by applying it to a simple case, the $4^4$ square  tiling seen in graph paper. 
The subgraph $H$ is shown in Fig.~\ref{FigGrid}.

\bigskip
\hrule
\bigskip
Figure~\ref{FigGrid} around here
\bigskip
\hrule
\bigskip

For this tiling it is  not hard to work out the CS directly, by 
drawing the contour lines, which are concentric squares centered at the origin.
Each new square plainly contains four more vertices than the last one,
so the tiling's CS satisfies the recurrence  $a(n+1) = a(n)+4$. 

Using the coloring book approach we draw ``trunks'' (in \textcolor{blue}{blue}) and ``branches'' (in \textcolor{OliveGreen}{green})
(see Fig.~\ref{FigGrid}). 
Each vertex is associated (via a blue or green edge) with a unique
vertex that is one step further away from $P$,
except for the vertices on the (blue) trunk,  where additional 
branches (green) have sprouted.\footnote{This example is slightly exceptional in that
the trunks and branches are not orthogonal to the contour lines.}
As four branches sprout at each distance $n$ from $P$
(the vertices marked with a green dot indicate the starts of the branches) 
the coordination sequence satisfies the recurrence  $a(n+1)= a(n)+4$.
The recurrence starts with $a(0)=1$, and so we have again shown that  $a(n)=1,4,8,12,16,20,\ldots$, as in \eqn{EqGrid}.


\section{The  Cairo tiling with respect to a tetravalent vertex}\label{SecTet}

At left in Figure~\ref{FigTet}, the vertices in the Cairo tiling in the vicinity of
a tetravalent vertex $P$ are marked by contour lines of constant distance to $P$, and counting the vertices  on each contour line confirms
that the initial terms of the CS are as shown in \eqn{EqGrid}.

\begin{figure} 
\centerline{\includegraphics[width=.9\textwidth]{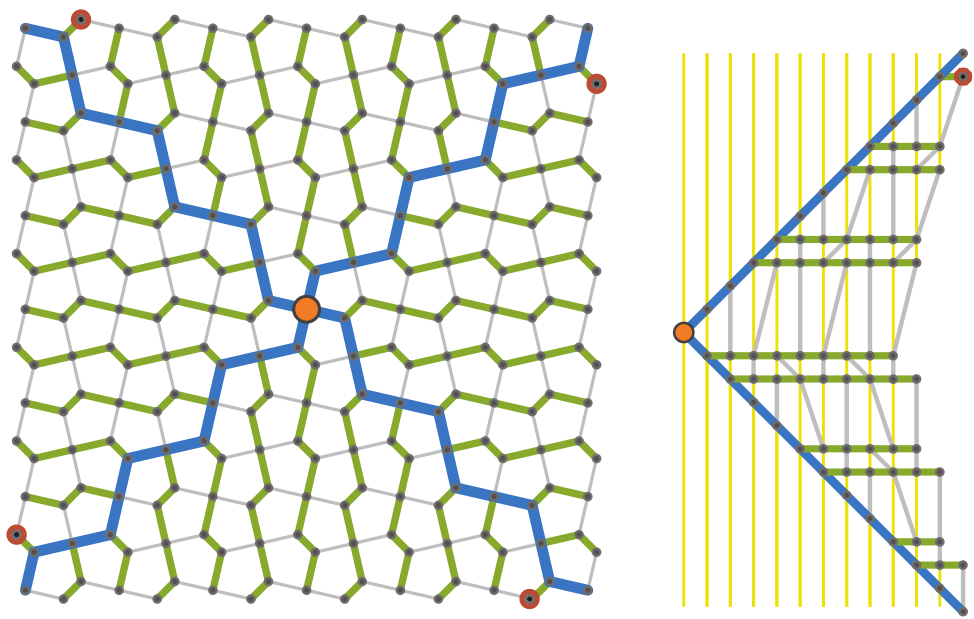}}
\caption{ At left, the subgraph $H$ for a tetravalent base
vertex in the Cairo tiling. There are four congruent sectors. At right, the  right-hand sector 
of $H$ has been redrawn so that the blue trunks and green branches 
are straight, with horizontal distances equal to the distance back to $P$ in $G$, in order to show that the (gray) edges not in $H$ do not reduce
the distance to the origin. Note that $H$ is transverse to the contour lines, shown in yellow.   At each distance $n+1$ from $P$,  all but four of the vertices (green) may be associated with the vertices at distance $n$, by tracing back one edge in $H$ --- hence 
$a(n+1)=a(n)+4$.
The four red dots on the periphery of the left figure
indicate the division into sectors, and help match up the left
and right figures. }
\label{FigTet2}\label{FigTet3}
\end{figure}

\bigskip
\hrule
\bigskip
Figure~\ref{FigTet2} around here
\bigskip
\hrule
\bigskip

The subgraph $H$ is shown in Fig.~\ref{FigTet2} (left). There are four trunks (blue) and
infinitely many branches (green), one branch originating at each trunk vertex. 
Figure~\ref{FigTet2} (right) shows one
sector of $H$ redrawn so that the trunks and branches are straight,
and so that we can see they satisfy Property (ii) --- that is, using
edges that are not part of $H$ (these edges are colored gray) does not provide any shorter
paths to the base vertex.

As each vertex on a trunk or branch is associated 
with another trunk or branch vertex  further out, and four new branches are introduced at each  distance $n>0$ from the origin, $a(n+1)= a(n)+4$; noting that $a(1)=4$, we have $a(n)=4n$, $n>0$, completing the proof of Theorem~\ref{ThTet}.

In the redrawn sector on the right
of Fig.~\ref{FigTet2}, all the points at a given distance from the base point
are colinear, and we see that there are exactly $n-1$ points in
the interior of the sector at distance $n$ from the base point.
Taking into account the four branch points, we have another proof
that $a(n) = 4(n-1) + 4 = 4n$ for $n>0$.

Note that there is a unique (green) branch originating at each (blue) trunk vertex.
If the trunk vertex is at distance congruent to $1$ or $2 \bmod{4}$ from the origin,
the branch turns to the left, otherwise it turns to the right.


\section{The Cairo tiling with respect to a trivalent vertex}\label{SecTri}

At right in Figure~\ref{FigTri} we show the contour lines in the vicinity of
a trivalent vertex $P$ in the Cairo tiling,  confirming that the initial terms of the CS 
are as shown in \eqn{EqA296368}.

\begin{figure} 
\centerline{\includegraphics[width=.84\textwidth]{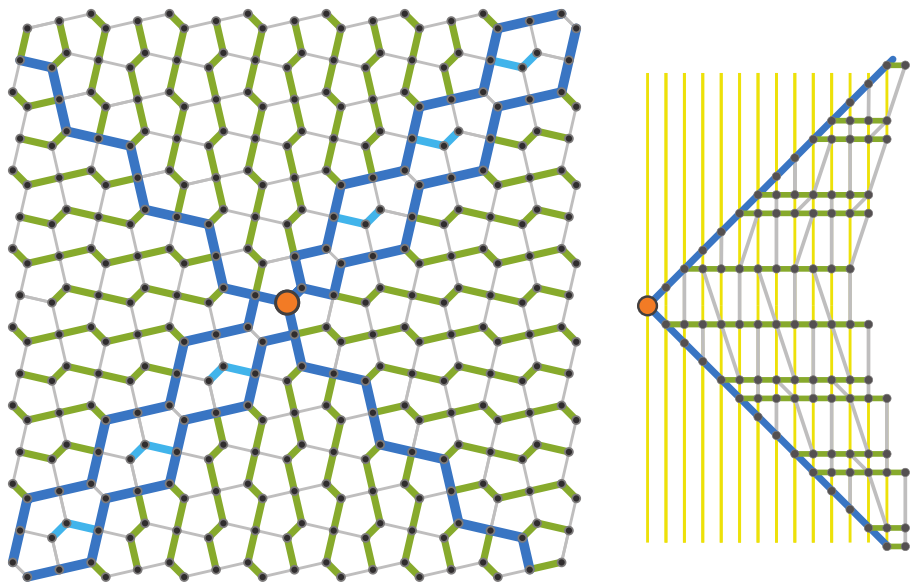}}
\caption{ At left, the subgraph $H$ (the blue, green, and  light blue 
edges) for the  trivalent case. Two of  the trunks 
have split down the middle, and moreover
have  sprouted ``twigs'' (light blue). 
In the right-hand illustration,
the right sector has been  redrawn so that the blue trunks and green branches 
are straight, in order to show that the (gray) edges not in $H$ do not reduce
the distances to the origin, and to make it easy to count the vertices  
at a given distance from the base point.} 
\label{FigTri2}\label{FigTri3}
\end{figure}

\begin{figure} 
\centerline{\includegraphics[width=.4\textwidth]{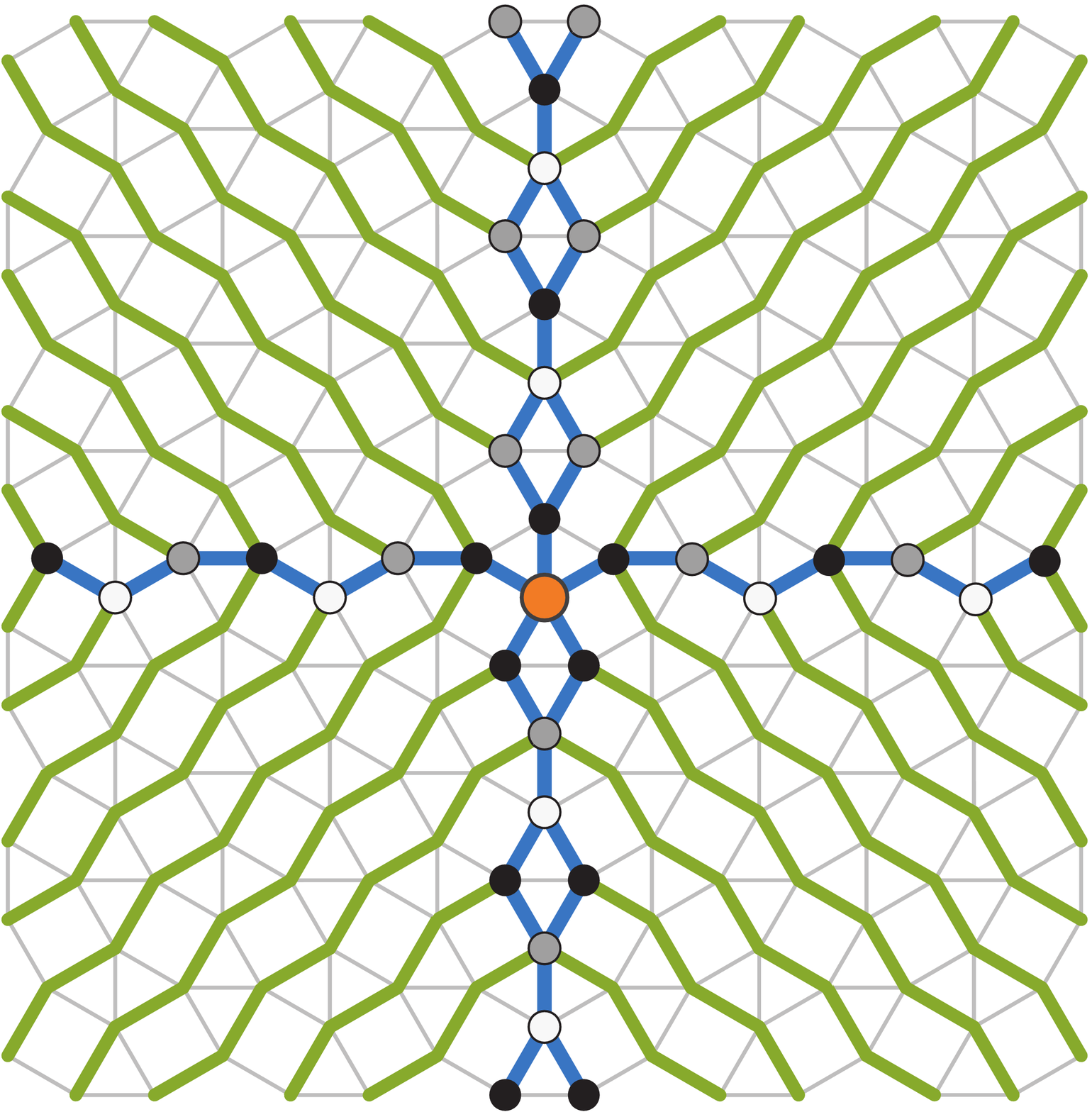}\hspace{.1\textwidth}\includegraphics[width=.4\textwidth]{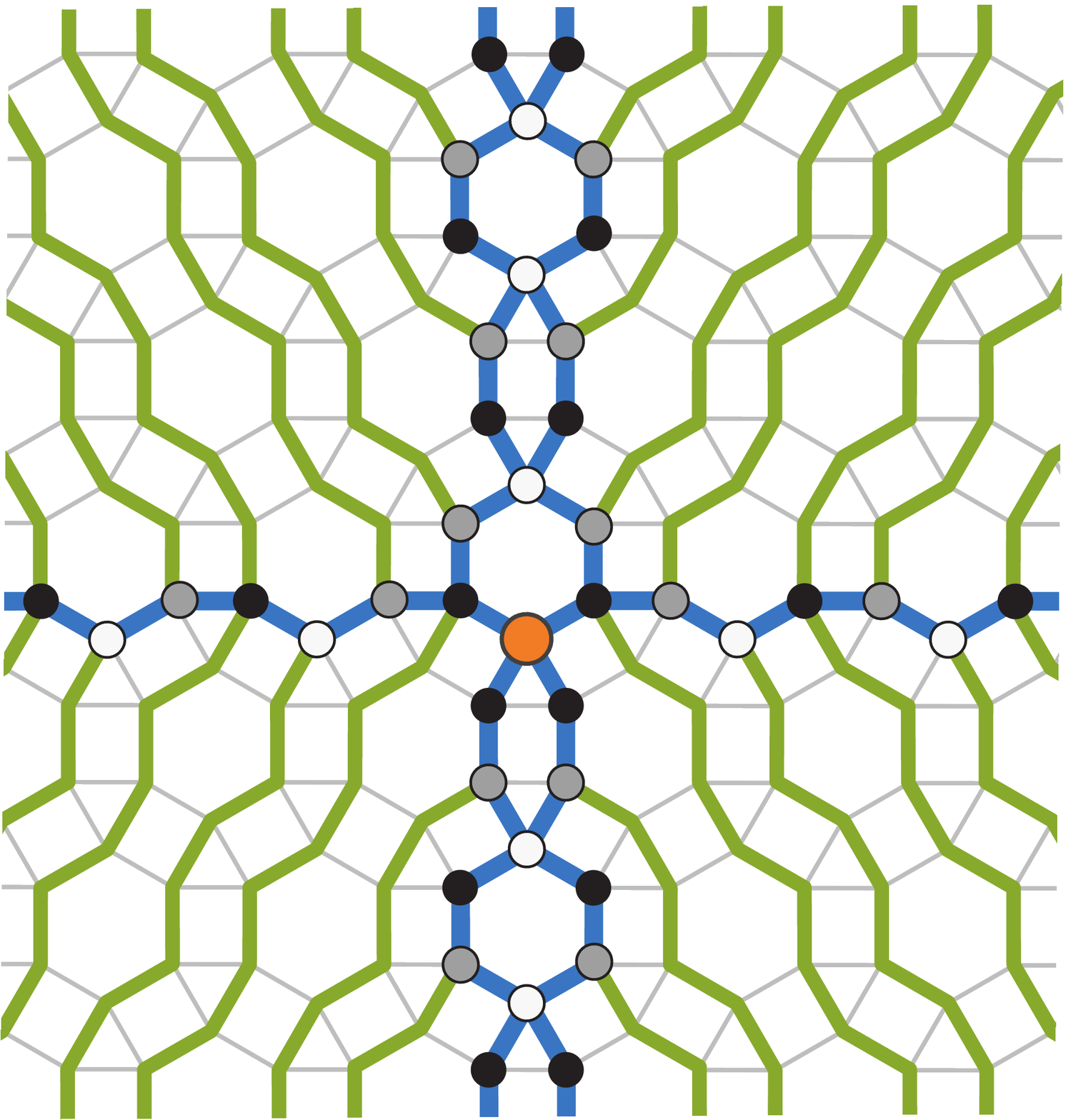}}
\caption{ Trunks and branches  $H$ for the $3^2 \ldot 4.3.4$ tiling, at left and for the $3.4.6.4$ tiling at right. In both drawings, the base vertex $P$ is at the center and vertices  on the blue trunks are colored by their distance to $P \mod 3$ --- $0$ (white), $1$ (black), or $2$ (gray).
}
\label{FigDualb}\label{Fig3464}
\end{figure}

Because the graph now has only mirror   symmetry,  the subgraph $H$ 
is necessarily less elegant than in the tetravalent case. 
The best choice for $H$ that we have found,
shown in Fig.~\ref{FigTri2},
now has six trunks, two of which sprout ``twigs'' shown in  light blue. 

\bigskip
\hrule
\bigskip
Figure~\ref{FigTri2} around here
\bigskip
\hrule
\bigskip

As in Fig.~\ref{FigTet2}, $H$ is naturally divided into four sectors (ignoring the twigs for now). The right and top sectors are congruent to each other, and the left and bottom sectors are essentially the same as any of the sectors in Fig.~\ref{FigTet2}, 
the main difference being that the base vertices  for 
the left and bottom sectors in Fig.~\ref{FigTri2} are now one edge away from the origin
instead of being at the origin as they were in Fig.~\ref{FigTet2}.

Although there is some variation from level to level, 
with two out of every three trunk nodes sprouting branches, 
and taking into account the periodic appearance of twigs,
we  obtain the recurrence $$a(n+4) =  a(n)+16,~ n\geq 3$$ by observing that exactly $16$ branches and twigs sprout in any four consecutive values
of $n$, and consequently all but $16$ vertices at distance $n+4$ may be traced back and associated with the vertices at distance $n$, along a trunk, or branch, or from a twig to a closer twig. 
Verifying the terms $a(0)$ through $a(6)$, we obtain $a(n)$ as in Theorem~\ref{ThTri}.

Again there is an alternative, more direct, way to see this. 
In Fig.~\ref{FigTri2} (right),
we see that the vertices at a given distance
from the base point are colinear, and for $n \ge 2$
there are $n-2$ interior vertices in the right (and top) sectors
if $n \equiv 0$ or $1 \pmod{4}$, or $n-1$ if $n \equiv 2$ or $3 \pmod{4}$.
In each of the left and bottom sectors there are $n-2$ interior vertices (for $n \ge 2$)
at distance $n$ from the base point. The numbers of trunk (or dark blue) vertices at distances $0,1,2,3,4,\ldots$ are $1,3,6,6,6,\ldots$, respectively,
and the numbers of twig (light blue) vertices at distances $4k$, $4k+1$, $4k+2$,
$4k+3$ (for $k \ge 1$) are $1,2,1,0$, respectively. 
Collecting these values,
we obtain equation \eqn{EqTrivalent}.


\section{The {$3^2 \ldot 4.3.4$} tiling}\label{SecDual}

The $3^2 \ldot 4.3.4$ tiling (Fig.~\ref{FigDual} in \S\ref{Sec1}) is  a uniform tiling: 
all vertices  are equivalent, and we can choose $P$ to be any  vertex. 
As in the previous section, the graph $G$ has only mirror
symmetry, and so again we have to accept
that our subgraph $H$ of trunks and branches will be less symmetrical that the graph in \S\ref{SecTet}.

In our choice for $H$, shown in Fig.~\ref{FigDualb},  there are two horizontal trunks, but on the vertical trunks some vertices  on the vertical branches have been split into ``burls''
(creating loops of length $4$).
With the orientation shown,  quadrants I and II are mirror images of each other,
as are quadrants III and IV. One can see immediately that any simple path in $H$ from a vertex back to $P$ is as short as any  other path back to $P$, and so Property (ii) is satisfied.

\bigskip
\hrule
\bigskip
Figure~\ref{Fig3464} around here
\bigskip
\hrule
\bigskip

To calculate the coordination sequence, we may  note that at any three 
consecutive distances $n+1$, $n+2$, $n+3$, $n>0$ from $P$, $16$ branches are introduced, and so all but $16$ vertices at distance $n+3$ may be associated with those at distance $n$. This gives us the recurrence $a(n+3)=a(n)+16, n>0$. Checking the  terms $a(0)=1, a(1)=5, a(2)=11, a(3)=16$, we obtain \eqn{EqDuala}, and so
complete the proof of Theorem~\ref{ThDual}.


\section{The {$3.4.6.4$} tiling}\label{Sec3464}

Our next example is the $3.4.6.4$ uniform tiling, and we will prove that this too has the same CS as the square grid
(this establishes a 2014 conjecture of Darah Chavey stated in \seqnum{A008574}).
\begin{thm}\label{Th3464}
The coordination sequence with respect to a vertex in
the $3.4.6.4$ tiling is given by $a(0)=1$, $a(n)=4n$ for $n \ge 1$.
\end{thm}

We show our choice of  trunks and branches at the right in Fig.~\ref{Fig3464}.
As in the left-hand figure, the vertical  trunks have  split, producing burls which
now are chains of hexagons. 
Again, quadrants I and II are mirror images of each other,
as are quadrants III and IV. We can see that Property (ii) is satisfied and that the pattern propagates. 

As a total of twelve branches are introduced at three
consecutive levels $n+1$, $n+2$, $n+3$, $n>1$, we have the 
recurrence $a(n+3)=a(n)+12$. Checking the terms up to $a(4) = 16$, we complete the proof.


\section{The {$4.8^2$} tiling}\label{Sec488}

\bigskip
\hrule
\bigskip
Figure~\ref{Fig488} around here
\bigskip
\hrule
\bigskip

\begin{thm}\label{Th488}
The CS (\seqnum{A008576}) 
with respect to a vertex in
the $4.8^2$ uniform tiling is given by $a(0)=1$
and thereafter $a(3k)=8k$, $a(3k+1)=8k+3$, $a(3k+2)=8k+5$ .
\end{thm}

\begin{proof}
The subgraph $H$, shown in at left in Fig.~\ref{Fig488} --- we again resort to using ``burls'' but the proof is the same as before: $H$ satisfies Property (ii) and may be propagated; for each $n>1$, at three consecutive distances a total of $24$ branches sprout. Consequently, for $n>1$, $a(n+3)=a(n)+24$. Verifying the terms through $a(4)$, we  complete the proof. \end{proof}


\bigskip
\hrule
\bigskip
Figure~\ref{Fig31212} around here
\bigskip
\hrule
\bigskip

\section{The {$3.12^2$} tiling}\label{Sec31212}

The $3.12^2$ uniform tiling is shown 
in Fig.~\ref{Fig31212}.
The CS with respect to any vertex begins
\beql{EqA250122}
1, 3, 4, 6, 8, 12, 14, 15, 18, 21, 22, 24, 28, 30, 30, 33, 38, 39, 38, 42, 48,  \dots \,.
\eeq

This is sequence \seqnum{A250122}, where there is a  conjectured formula from 2014
due to Joseph Myers which we can now prove is correct.

\begin{thm}\label{Th31212}
The coordination sequence for the $3.12^2$ tiling is given by
$a(0)=1$, $a(1)=3$, $a(2)=4$, and thereafter $a(4k)=10k-2$, 
$a(4k+1)=9k+3$, $a(4k+2)=8k+6$, and
$a(4k+3)=9k+6$.
\end{thm} 

\begin{proof}
The subgraph $H$ is shown is shown 
in Fig.~\ref{Fig31212}. 
There are double trunks along the $x$-axis, and a single trunk with burls
along the $y$-axis.
The second quadrant is a mirror image of the first,
and the third quadrant is a mirror image of the fourth, so we need only analyze 
the first and fourth quadrants. (The figure is not symmetric
about the $x$-axis.) To simplify the discussion we assume the positive $y$-axis
points North and the positive $x$-axis points East.

In the first quadrant there are infinitely many parallel branches (green)
consisting of infinite paths directed to the North-East, together with twigs (olive green)
of length $2$ originating at certain branch vertices. The first quadrant is bisected by 
one of the branches, let us call it the special branch, indicated by a thick red line, which roughly follows the diagonal $y=x$.
To the West of the special branch, the twigs are all directed to the North-West, while to the East of the special branch the twigs are directed to the East.

The situation in the fourth quadrant is similar. 
There are infinitely many parallel branches (green)
consisting of infinite paths directed to the South-East, together with twigs (olive green)
of length $2$ originating at certain branch vertices. The fourth quadrant is bisected by 
one of the branches, the special branch, indicated by a thick red line, which roughly follows the diagonal $y=-x$.
To the West of the special branch, the twigs are all directed to the South-West, while to the East of the special branch the twigs are again directed to the East.

\bigskip
\hrule
\bigskip
Table~\ref{Tab31212} around here
\bigskip
\hrule
\bigskip

Table \ref{Tab31212} shows the numbers of vertices of the various types at distance $n$ from the central vertex. The counts depend on the value of $n$ modulo $8$,
as indicated by the columns of the table. We assume $n \ge 3$.
The rows of the table indicate the different types of vertex.
Row (i) refers to the East-West trunks, row (ii) to the North-South trunk,
excluding the vertices already counted in (i). 
Row (iii) counts vertices in the first quadrant that are on and to the East of the
special branch, (iv) those to the West of that branch, and (v) counts the twigs.
Rows (vi), (vii), (viii) are the analogous counts for the fourth quadrant.The last row of the table gives the grand total, formed by adding the entries in rows (i) and (ii),
plus twice (to account for the second and third quadrants) the sum of rows (iii) though (viii). 

Examination of the last row shows that the total actually depends only on the value of $n$ modulo $4$, rather than $8$. For example,
the assertions that $a(8k)=20k-2$ and $a(8k+4)=20k+8$
can be combined into $a(4k)=10k-2$.
Similarly for the other six cases. This completes
the proof of the theorem.
\end{proof}


\section{The  {$3^4 \ldot 6$} tiling}\label{Sec33336}

The $3^4 \ldot 6$ uniform tiling is shown in Fig.~\ref{Fig346a}.
This is also known as the snub $\{6,3\}$ tiling, and has symmetry group $632$.
The CS with respect to any vertex begins
\beql{EqA250120}
1, 5, 9, 15, 19, 24, 29, 33, 39, 43, 48, 53, 57, 63, 67, 72, 77, 81, 87, 91, \dots \,.
\eeq
This is sequence \seqnum{A250120}, where there is a  conjectured formula from 2014
due to several people, which we can now prove is correct.

\bigskip
\hrule
\bigskip
Figure~\ref{Fig346a} around here
\bigskip
\hrule
\bigskip

\begin{thm}\label{Th33336}
The coordination sequence for the $3^4 \ldot 6$ tiling is given by
$a(0)=1$, $a(1) = 5$, $a(2) = 9$, $a(3)=15$, $a(4)=19$, $a(5) = 24$ 
and  for $n \ge 3$, $a(n+5)=a(n)+24$. 
\end{thm} 
The recurrence can also be written as $a(5k+r)=24k+a(r)$ for $k \ge 0$, $0 \le r \le 4$.

In all our previous examples, we were able to use the coloring
book method to find trunks and branches by simply drawing on
a picture of the tiling. This enabled us to calculate the coordination
sequence directly, bypassing any complexities in the contour lines.
Compare, for example, the contour lines in Fig.~\ref{FigTet}
with the trunks and branches of Figures~\ref{FigTet2} and~\ref{FigTri2}, 
which made the CS for the Cairo tiling a simple calculation.
We were happy that this strategy  has worked 
in every case up to this point.

Of course, to use the coloring book method,
we must also confirm, in some inductive manner,  
that the trunks and branches can be continued indefinitely  
and satisfy Property (ii). For the Cairo tiling, we redrew two of the sectors
so as to make it clear that there were patterns
of local structure that propagated outwards, maintaining a valid trunks
and branches structure that made it easy to compute the CS.
We were not as explicit in the next few examples, leaving to the reader  the
easy verification that the trunks and branches structure  propagated.   
In the present section and the next, however, the structures
are more complicated and so we must be more cautious. 

We could have always  used the contour lines directly, at least in principle ---  
defining the contour lines recursively,  with production rules that are
constrained by local conditions, as in Goodman-Strauss (2009). For tilings with a co-compact symmetry, such as the ones in this paper, contour lines can be used to conjecture, 
calculate, and give proofs for the recurrences
satisfied by the coordination sequences.

The contour lines for the $3^4 \ldot 6$ tiling are shown in  Fig.~\ref{Fig33336Level}.
From this we can verify that, at least initially,  the number of vertices at distance $n>1$
from the base vertex increases by exactly $24$ when $n$ is increased by $5$.
Although there appear to be a lot of individual cases to consider,
the contours are highly ordered,  and the task is 
primarily one of cataloging the local structures.

\bigskip
\hrule
\bigskip
Figure~\ref{Fig33336Level} around here
\bigskip
\hrule
\bigskip

Outside of an initial  region about $P$, there are really only four 
kinds of vertex neighborhoods to 
consider: vertices where the contour lines are straight, where they zig-zag, and where
there is a transition from zig-zag-to-straight and {\it vice versa}.
In each case, the underlying structure is straightforward and it is
easy to see that it propagates from level to level, and also that the structure we see at the center --- twelve sectors of alternating types, the six sectors of each type being fundamentally the same --- must continue across the entire infinite tiling. 

However, there is a further  complication:
the boundaries between straight-to-zig-zag,
drawn in black in Fig.~\ref{Fig33336Level}, have a  period-$3$ pattern, and so
the structure of the segments crossing each sector recurs with period $15$.

We now use  a  trunks and branches approach to simplify this structure, verifying locally 
that this give the CS, but without needing to check that there is a detailed
match with the contour lines.  

\bigskip
\hrule
\bigskip
Figure~\ref{Fig33336Coloring} around here
\bigskip
\hrule
\bigskip

The trunks and branches structure  is shown on the left in Fig.~\ref{Fig33336Coloring}.
As in Fig.~\ref{Fig3464}, the base vertex $P$ is at the center and vertices  on the blue trunks are colored by their distance to $P \mod 3$ --- $0$ (white), $1$ (black), or $2$ (gray).
The six trunks each produce a pair of branches and a twig at every third level, but shifted in pairs: all together the trunks sprout four branches and four twigs at each level. 
On the right in Fig.~\ref{Fig33336Coloring}
are shown four patterns (1)-(4) that are self-perpetuating from level to level. 
If the boundary of a disk can be described by this atlas of patterns, then it is inside a larger disk with the same property --- though one must carefully check that every transition from zig-zag-to-straight  is of this precise form.  By induction, the alleged contour lines in the pattern are nested simple closed curves and so can only lie on actual contour lines; thus  the trunks and branches satisfy Property (ii) and  correctly give the CS. With the exception of the red and yellow vertices  in pattern (1),  each vertex at level $n+5$ is naturally matched with a vertex at level $n$, just by following a branch backwards. Notice that the branches have matching twigs at every fifth level. In pattern (1), we examine the mismatched vertices: the red vertices  have no match five levels earlier,  and the yellow vertex has no match five levels forward. Since pattern (1) repeats at every third level, we have considered all possible cases. At each level,  two of the trunks show each of the three cases in pattern (1),  so the number of vertices  satisfies   $a(n+5) = 24+a(n)$ for $n\geq 3$, a net increase of six vertices  per trunk. 
After explicitly verifying the initial terms of the coordination sequence, we have completed the proof of Theorem~\ref{Th33336}.


\section{The snub-632 tiling}\label{SecSnub}
 
\bigskip
\hrule
\bigskip
Figure~\ref{Figsnub1} around here
\bigskip
\hrule
\bigskip

We began this article by  analyzing the Cairo tiling,
which is the dual to the $3^2 \ldot 4.3.4$ uniform tiling.
Our final example is the snub-$632$ tiling (Fig.~\ref{Figsnub1}), which is
the dual to the $3^4 \ldot 6$ uniform tiling of the previous section.
Like the Cairo tiling, this is a beautiful tiling; fundamentally it is the dual of
the snub versions of both the $6^3$ and $3^6$ tilings.
It has several names including: \begin{itemize}
\item the $6$-fold pentille  tiling 
(Conway et al., 2008, p.~288),
\item the fsz-d net 
(O'Keeffe et al. 2008),
\item and  the dual of the $3^4 \ldot 6$ tiling 
(Gr\"{u}nbaum  \& Shephard, 1987, pp.~63, 96, 480 (Fig.~$P_5$-16)).
\end{itemize}
We will refer to it simply as the snub-$632$ tiling.
There is only one shape of tile, an elongated pentagon.\footnote{Again,
as long as the underlying topology of the  graph is preserved, small variations in the ratios of the sides do not affect the coordination sequences.}
There are three types of vertices:  
(a) hexavalent vertices, with six-fold rotational 
symmetry, where six  long edges of the pentagons meet;
(b) trivalent vertices,  with three-fold rotational symmetry,  where three short edges meet; 
and
(c) trivalent vertices, with no symmetry,  where two short edges and one long edge meet.
The coordination sequences (\seqnum{A298016}, \seqnum{A298015}, 
\seqnum{A298014}) for these three types of vertices begin
\begin{align}\label{Eqsnub1}
(a)~~ & 1, 6, 12, 12, 24, 36, 24, 42, 60, 36, 60, 84, 48, 78, 108, 60, 96, 132, 72, \ldots \,, \nonumber \\
(b)~~ & 1, 3, 6, 15, 24, 18, 33, 48, 30, 51, 72, 42, 69, 96, 54, 87, 120, 66, 105,  \ldots \,, \nonumber \\
(c)~~ & 1, 3, 9, 15, 18, 27, 37, 37, 44, 57, 54, 61, 77, 71, 78, 97, 88, 95, 117, 105,  \ldots 
\end{align}

\begin{thm}\label{Thsnub}
The CS for the three types of vertices  in the snub-$632$ tiling are given by: \\
(a) $a(3k)=12k,~ a(3k+1)=18k+6,~ a(3k+2)=24k+12$ for $k\ge 1$, \\
(b)  $a(3k)=18k-3,~ a(3k+1)=24k,~ a(3k+2)=12k+6$ for $k\ge 2$, \\
(c)  $a(3k)=20k-3,~ a(3k+1)=17k+3, ~a(3k+2)=17k+10$ for $k\ge 2$, \\
with initial values as in \eqn{Eqsnub1}. In each case
we have $a(n)=2a(n-3)-a(n-6)$ for $n \ge 8$.
\end{thm}

\bigskip
\hrule
\bigskip
Figure~\ref{Fig632SnubCoordination} around here
\bigskip
\hrule
\bigskip

We could prove Theorem~\ref{Thsnub} by working directly with the contour lines, shown  in Figure~\ref{Fig632SnubCoordination}, and
establishing that there is a recursive --- although unwieldy ---  structure for each of the three
types of vertices.  However,  the coloring book approach will enable us to give  a 
uniform treatment.

\bigskip
\hrule
\bigskip
Figure~\ref{Fig:snubColoring2} around here
\bigskip
\hrule
\bigskip

Figure~\ref{Fig:snubColoring2} shows the trunks and branches structures 
for the three types of vertices, which we
will continue to refer to as cases (a), (b), and (c) for a $6$-fold, 
$3$-fold, and asymmetric base vertex, respectively. Each figure is divided into $60$-degree
sectors, some of which are separated by channels
(pairs of parallel branches distance $3$ apart).
The individual sectors are
all the same, although they differ in how far the apex of the sector (yellow)
is from the base vertex (red). In case (a) all six sectors have apex at the base vertex,
in case (b) the base vertex is two steps away from the apex, in two different ways,
and in case (c) the base vertex is  one, two or three steps away from
the apex, in six different ways (now each sector has a different displacement from the base vertex). 

\bigskip
\hrule
\bigskip
Figure~\ref{Fig:snubColoring1} around here
\bigskip
\hrule
\bigskip

Figure~\ref{Fig:snubColoring1} shows a sector (internally they are all the same)
with a channel next to it. The sector is bounded by two trunks (blue) and 
has a pattern of branches (red) and twigs (green). The figure also shows
the contour lines defined by the distances to the apex (yellow).
Contours at distances congruent to $0$, $1$, and $2$ mod $3$ 
from the apex are colored light gray, dark gray, and black, respectively.
A crucial step in the analysis is 
that these contour lines are still the contour
lines with respect to the base vertex, even when the base vertex 
is at one of the other eight possible vertices --- there are no shortcuts to the base vertex.
Of course the distances from the contours
to the base vertex get increased by one, two, or three steps
when the base vertex is moved. We return to this point
in Fig.~\ref{Fig:632snubcontourdetail}.

In the sector, the vertices  at distance $n+3$ from the apex are in 
 one-to-one correspondence with the vertices  at  distance $n$, with the exception 
 (and this is the other crucial step)
 of
 the nine vertices shown in red on the right of Fig.~\ref{Fig:snubColoring1}.
There are two unmatched vertices  at distance $0 \bmod{3}$ from the apex  (on the light gray curve), three at distance $1 \bmod{3}$ (on the dark gray curve)  and four at distance $2 \bmod{3}$ (on the black curve). As the structure along a trunk repeats with period $3$, this includes all the possibilities.
The vertices  that are not in any sector are in exact correspondence with those three
steps  further away.  

The detailed book-keeping for case (c) is as follows; (a) and (b) are simpler and are
discussed below.

In case (c), take  $n$ sufficiently large, which turns out to mean  $n > 5$, 
that being the distance from $P$ to  where the regular structure starts to be self-perpetuating. 
There is one sector with apex at distance $3$; so from P, there are, for 

\medskip
\begin{tabular}{ll}
$n \equiv 0 \mod 3$, & two unmatched vertices at distance $n+3$, \\
$n \equiv 1 \mod 3$, & three unmatched at $n+3$,\\
$n \equiv 2 \mod 3$, & four unmatched at $n+3$.
\end{tabular}
\medskip

There are two sectors with apex at distance $2$, so from P there are, for 

\medskip
\begin{tabular}{ll}
$n \equiv 2 \mod 3$,  twice two unmatched vertices at distance $n+3$, \\
$n \equiv 0 \mod 3$,  twice three unmatched at $n+3$, \\
$n \equiv 1 \mod 3$,  twice four unmatched at $n+3$.
\end{tabular}
\medskip

There are three sectors with apex at distance $1$, so from P there are, for  

\medskip
\begin{tabular}{ll}
$n \equiv 1 \mod 3$, three times  two unmatched vertices at distance  $n+3$, \\
$n \equiv 2 \mod 3$,  three times    three unmatched at $n+3$, \\
$n \equiv 0 \mod 3$, three times     four unmatched at $n+3$ 
\end{tabular}
\medskip

In summary,  for case (c) there are three sectors with apex at distance $1$ from $P$, two at distance $2$ and one at distance $3$. The vertices in the channels do not contribute to the recursion (since the vertices at distance $n+3$ are in one-to-one correspondence with those at distance $n$, for $n>2$). For case (c) we therefore have  the recurrence, for $k>1$, 

\medskip
\begin{tabular}{ll} 
$a(3k+6) = 3\cdot 4 + 2\cdot 3+ 1\cdot 2+ a(3k+3)$,  \\
$a(3k+7) = 3\cdot 2 + 2\cdot 4+ 1\cdot 3+ a(3k+4)$, \\
$a(3k+8) = 3\cdot 3 + 2\cdot 2+ 1\cdot 4+ a(3k+5)$. 
\end{tabular}
\medskip

In case (a), there are six sectors, each with apex at distance $0$ from $P$. The recurrence is therefore, for $k>1$, 
$a(3k+3) = 6\cdot 2 + a(3k)$,  
$a(3k+4) = 6\cdot 3 + a(3k+1)$,
$a(3k+5) = 6\cdot 4 + a(3k+2)$.

For case (b), there are six sectors, each with apex at distance $2$ from $P$, 
and the recurrence is, for $k>1$, 
$a(3k+5) = 6\cdot 2 + a(3k+2)$,  
$a(3k+6) = 6\cdot 3 + a(3k+3)$, 
$a(3k+7) = 6\cdot 4 + a(3k+4)$.

After verifying the initial terms by hand, we recover the sequences stated in the theorem. 

However, as mentioned above, we must still check 
that even when the base vertex $P$ is not at the apex,
the contour lines measured from $P$ really do cross the sectors and their trunks and branches structures as shown in Figure~\ref{Fig:snubColoring1}.
This check is carried out in Figure~\ref{Fig:632snubcontourdetail}. 
The heavy lines indicate where the regular structure --- the recurrence --- begins. 

\bigskip
\hrule
\bigskip
Figure~\ref{Fig:632snubcontourdetail} around here
\bigskip
\hrule
\bigskip


\section{Cayley diagrams and growth series}\label{SecCD}
\begin{figure} 
\centerline{\includegraphics[width=.4\textwidth]{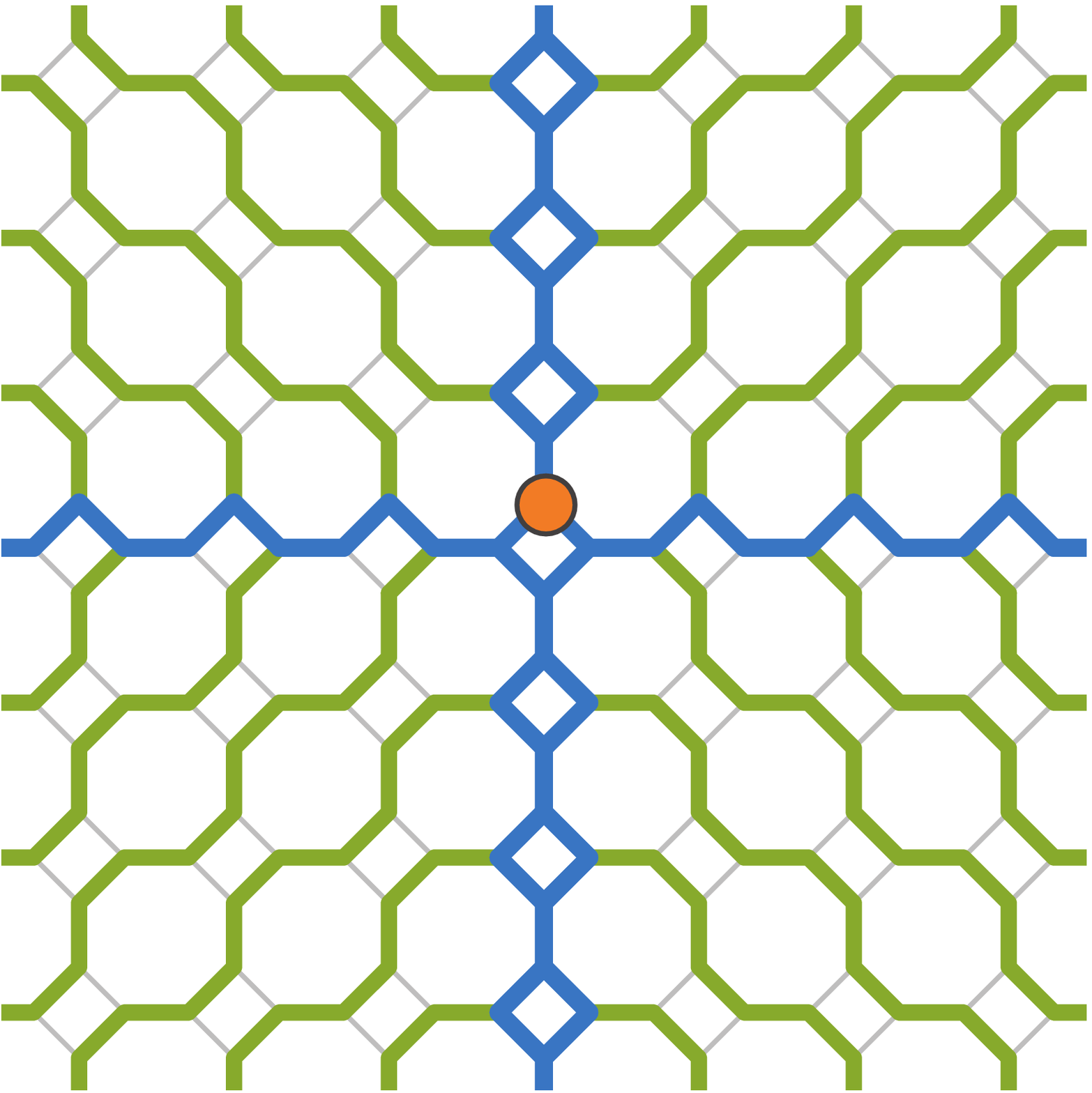}}
\caption{Trunks and branches $H$ for the $4.8^2$ tiling, 
with $P$ at the center and trunks and branches $H$ which clearly propagate and satisfy Property (ii). }
\label{Fig488} 
\end{figure}

\begin{figure} 
\centerline{
\includegraphics[width=.4\textwidth]{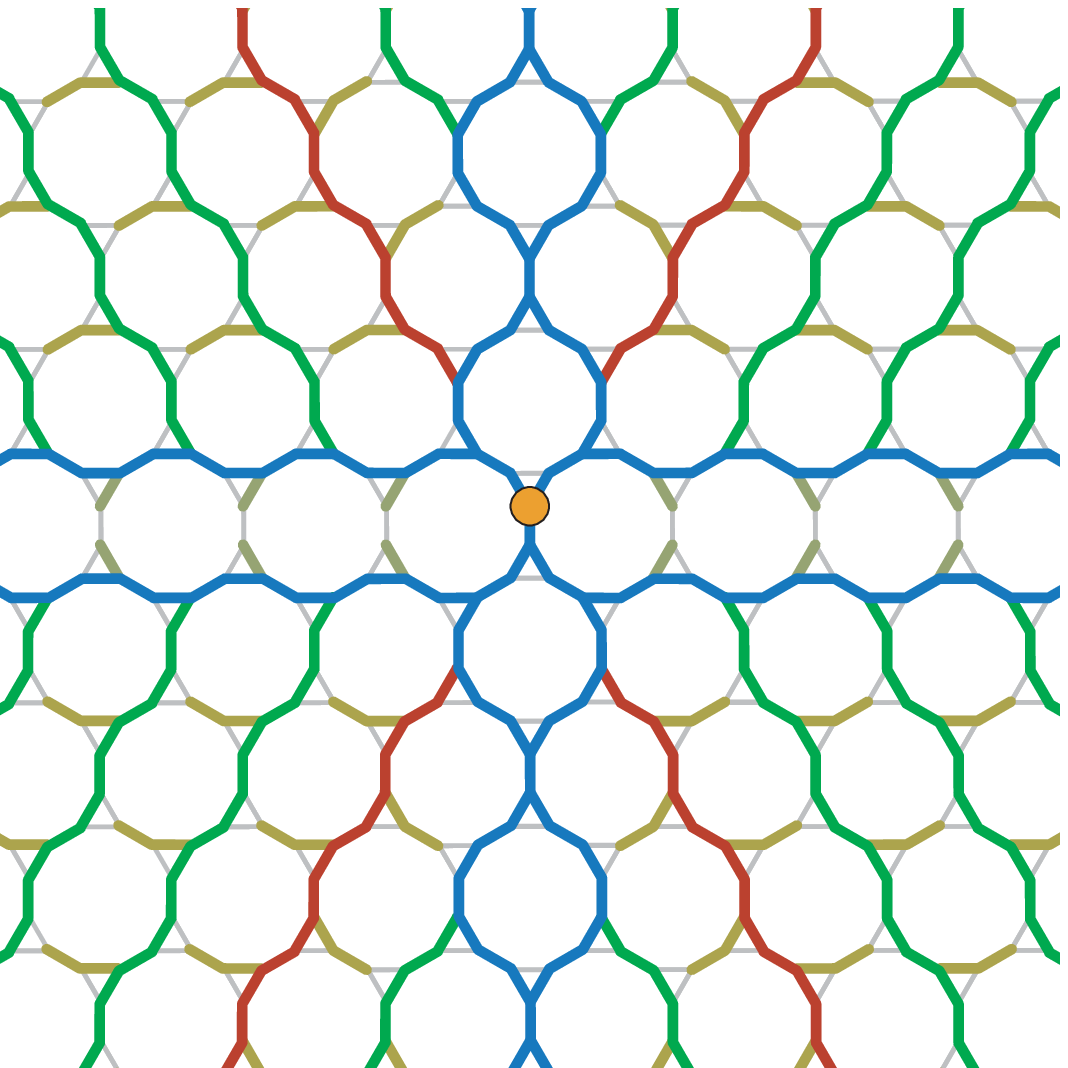}}
\caption{Trunks and branches $H$ for the  $3.12^2$ tiling, with $P$ at the center and trunks and branches $H$ which clearly propagate and satisfy Property (ii). }
\label{Fig31212}
\end{figure}

\begin{figure} 
\centerline{\includegraphics[width=.9\textwidth]{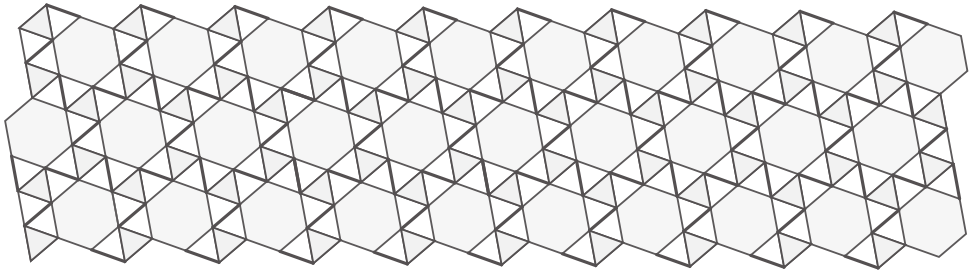}}
\caption{A portion of the $3^4 \ldot 6$ uniform tiling.}
\label{Fig346a}
\end{figure}

\begin{figure} 
\centerline{\includegraphics[width=.6\textwidth]{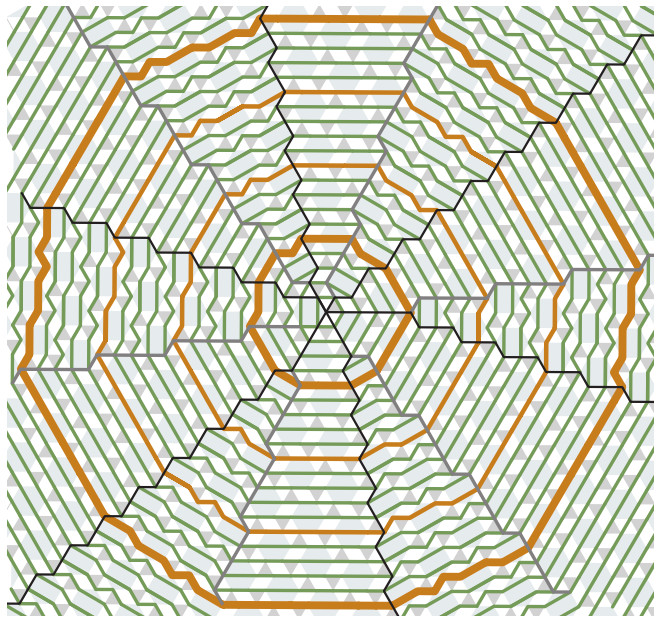}}
\caption{ The contour lines of the  $3^4 \ldot 6$ tiling; their lengths form sequence \seqnum{A250120}. The contour lines are partitioned into twelve sectors, alternatively straight and zig-zag. The boundaries   zig-zag-to-straight (reading clockwise), drawn in black, have period $3$. The boundaries straight-to-zig-zag, drawn in gray, have period $5$ --- in each sector the patterns at the ends of a contour segment have period $15$. A pair of contours $15$
steps  apart are emphasized in the figure by heavy brown lines. With patience, after another $15$ steps   
the regular structure becomes more apparent.
Even without doing that, with care  one can verify  that from each level $n>1$ to level $n+5$, the number of vertices  increases by $24$. In Fig.~\ref{Fig33336Coloring}, we use a trunks
and branches structure to make this more transparent.}
\label{Fig33336Level}
\end{figure}

\begin{figure} 
\centerline{\includegraphics[width=\textwidth]{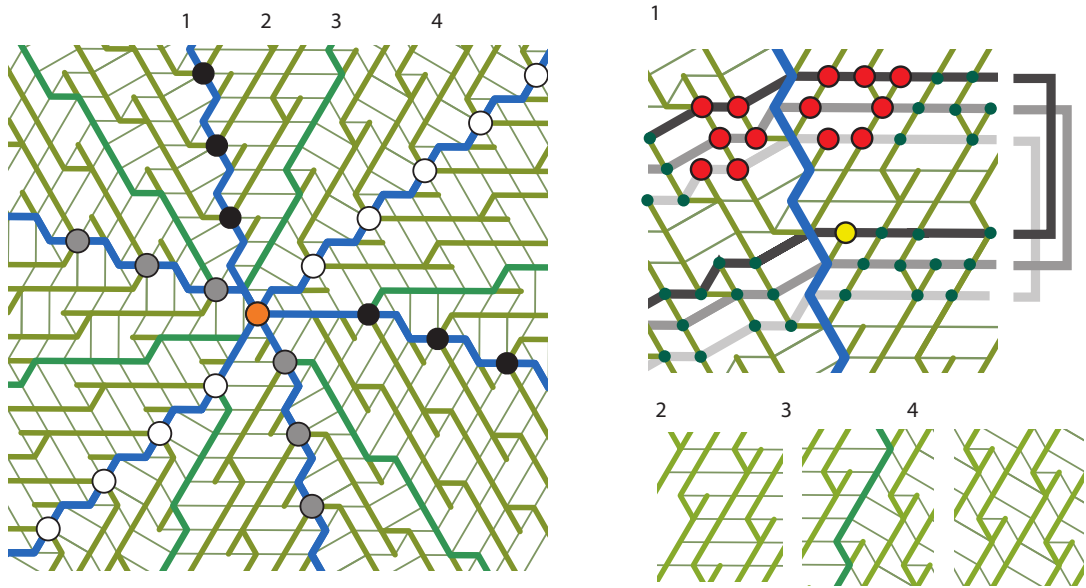}}
\caption{A trunks and branches structure for the $3^4 \ldot 6$ tiling (left).
The four patterns (1)-(4) that self-perpetuate from level to level (right).
In pattern (1), the red vertices at each distance $n+3$ are unmatched at distance $n$, and the yellow vertex at distance $n$ is unmatched at distance $n+3$. For sufficiently large $n$ ($>4$) there are a total of $12$ more vertices at distance $n+3$ than at $n$. }
\label{Fig33336Coloring}
\end{figure}

\begin{figure}
\centerline{\includegraphics[width=.8\textwidth]{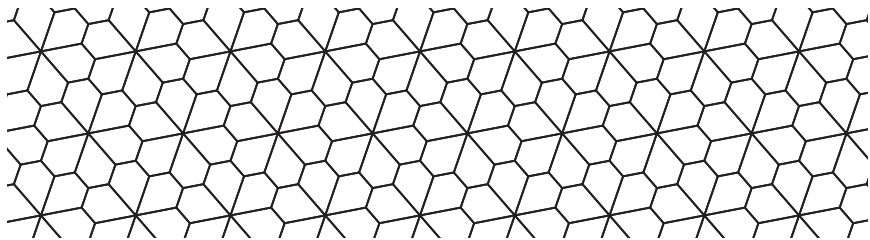}}
\caption{The snub-$632$ tiling. }
\label{Figsnub1}
\end{figure}

\begin{figure}
\centerline{\includegraphics[width=\textwidth]{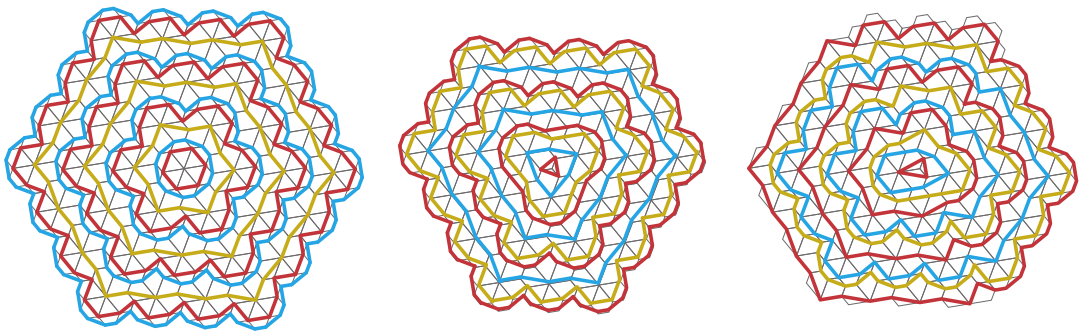}}
\caption{Contour lines in the snub-$632$, centered on a six-fold, three-fold and asymmetric vertex; they are clearly structured and we could compute with them directly. But the  
coloring book approach simplifies and unifies the  count. } 
\label{Fig632SnubCoordination}
\end{figure}

\begin{figure}
\centerline{\includegraphics[width=.8\textwidth]{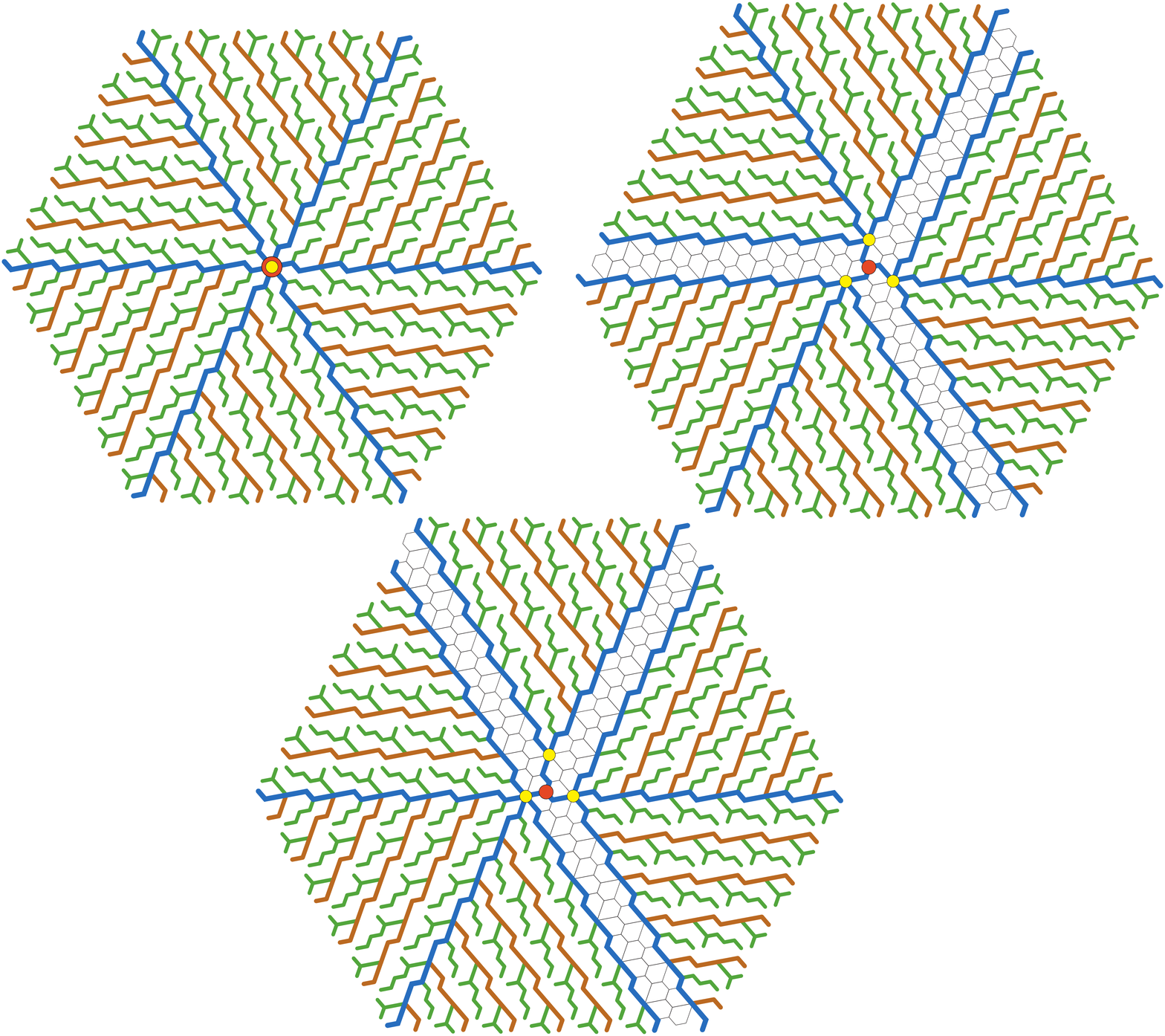}}
\caption{ Trunks and branches structures for the three types of base vertices in the
the snub-$632$ tiling: (a) (top left), 6-fold base vertex, (b) (top right), 3-fold, (c) (bottom) asymmetric. 
Each figure is divided into $60$-degree
sectors, some of which are separated by channels. 
The base vertex $P$ is shown in red and the vertices at the apices of the sectors in yellow.
}
\label{Fig:snubColoring2}
\end{figure}

\begin{figure}
\centerline{\includegraphics[width=.7\textwidth]{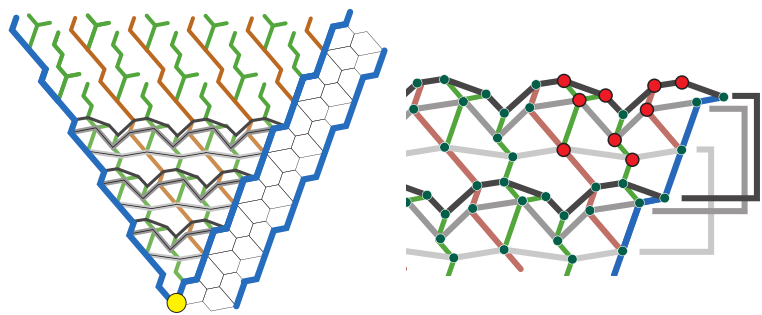}}
\caption{On the left, a sector and a channel,
showing trunks (blue), branches (red), twigs (green), and contours
(light grey, dark gray, and black) defined by
distances to the apex (yellow).
 In the sector, the vertices  at distance $n+3$ from the apex are in 
 one-to-one correspondence with the vertices  at  distance $n$, with the exception of those shown in red at right.
 Within each sector, then, for $n$ sufficiently large (that is, beyond the tip of the sector),  we obtain a recursion for the number of vertices at distance $n+3$, depending on $n$ modulo $3$. For $n\equiv 0$ there are two more vertices than at distance $n$, for $n\equiv 1$ there are four more, and for $n\equiv 2$ there are three more. In a channel  there are the same number of vertices at distance $n+3$ as for $n$. In the text we complete the calculation across the whole tiling.}
\label{Fig:snubColoring1}
\end{figure}

\begin{figure}
\centerline{\includegraphics[width=.92\textwidth]{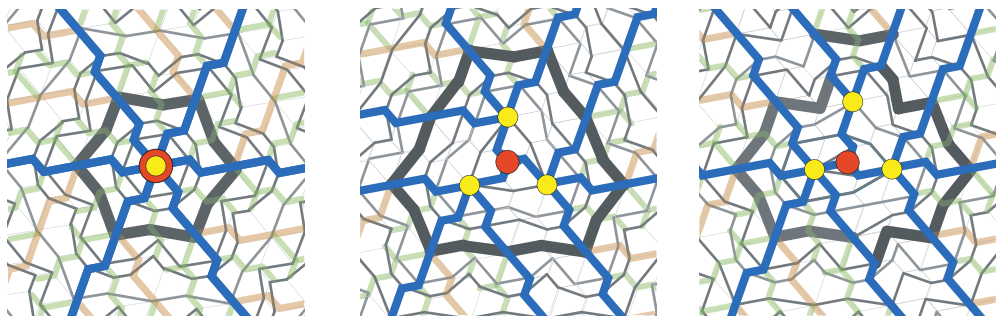}}
\caption{The contour lines near $P$ in each of the cases (a), (b) and (c).}
\label{Fig:632snubcontourdetail}
\end{figure}

If the symmetry group of the tiling acts transitively on the vertices,
and the subgroup fixing a vertex is the trivial group,  it is possible to represent the graph of the tiling as a Cayley diagram 
for an appropriate presentation of the group. For example,
consider the $4.8^2$  tiling shown in Fig.~\ref{Fig488}.
Every vertex is trivalent, with one edge (denoted by `$o$', say) separating
two adjacent octagons, and a pair of edges (`$\ell$' and `$r$') 
that go either to the left or the right around the adjacent square.
Let $R$ mean `move along the edge $o$',
and let $S$ mean `move along the edge $\ell$'.
Applying $R$ twice returns one to the start, which we indicate
by saying that $R^2=1$. Similarly, $S^4=1$ (going around the square)
and $(RS)^4=1$ (going around an octagon).
The group $G$ generated by $R$ and  $S$ subject
to these relations, specified by the presentation
\beql{EqPres}
G ~=~ \langle R, S \mid R^2 = S^4 = (RS)^4 = 1 \rangle\,,
\eeq 
assigns a unique label to each vertex in the graph.  
In technical terms, this gives an identification of the graph
of the tiling with the Cayley diagram of the group defined by this
presentation (Johnson 1976).

The `length' of a group element is the minimal number of
generators needed to represent it, which is also its distance from
the identity element in the Cayley diagram.
The `growth function' for the group specifies the number of
elements of each length $n$, and so the
growth function is precisely the coordination sequence 
for the tiling. For the graphs of uniform tilings,
the Knuth-Bendix algorithm
(Knuth \& Bendix, 1989),  Epstein et al., 1991)  
can be used to solve the 
word problem\footnote{Although in general this problem is insoluble.} and
determine the growth function. This algorithm  is implemented
in the GrowthFunction command in the computer algebra system
Magma (Bosma et al. 1997).
When applied to the presentation~\eqn{EqPres}, for
example, Magma returns the generating function
$$
\sum_{n=0}^{\infty} a(n) x^n ~=~ \frac{(1+x)^2}{(1-x)^2}\,,
$$
which is equivalent to the formulas given in Theorem~\ref{Th488}.

\bigskip
\hrule
\bigskip
Table~\ref{TabCox} around here
\bigskip
\hrule
\bigskip

Table~\ref{TabCox} lists presentations for all eleven uniform 
two-dimensional tilings. We have verified that the resulting growth functions 
coincide with the coordination sequences given in earlier sections
and in the OEIS.

Presentations for the $17$ planar crystallographic groups
can be found in Table~3 of Coxeter \& Moser (1984), 
and  Shutov (2003) used them
to compute coordination sequences for the corresponding  {\em directed}
Cayley graphs. 
Eon (2018) gives presentations for the eleven uniform
tilings and points out the connection
between the Cayley diagrams and the graphs of the tilings.
However, neither article uses these presentations to explicitly
compute the coordination sequences for the tilings.

It is worth pointing out that the growth series approach to finding
coordination sequences used in this section only applies when the graph of the tiling 
coincides with the Cayley graph of some group, whereas our
coloring book approach can potentially be applied to any periodic tiling.


\ack{Acknowledgments}
We thank  Davide M. Proserpio for telling us about the {\em RCSR} and
{\em ToposPro} websites, and for drawing our attention to 
the articles Eon (2018) and Shutov (2013).
Davide M. Proserpio has also helped by using {\em ToposPro}
to compute coordination sequences for many tilings,
not mentioned in this article, which are now
included in the OEIS.

After seeing a preliminary version of this
article,  Jean-Guillaume Eon commented that
the coloring book method might be viewed as complementary
to the algebraic method introduced in Eon (2002).
The labeled quotient graph of the net or tiling
from that article could help to find the subgraph $H$ needed for
our approach, and conversely, using $H$ instead of the quotient graph
may simplify finding the the coordination sequence.

It will be interesting to see if a combination of our methods
will lead to proofs of the conjectured formulas for some of the more
complicated tilings, such as those of the $20$ $2$-uniform tilings.
A list of these $2$-uniform tilings and their
coordination sequences (and conjectured formulas)
may be found in entry \seqnum{A301724} in the OEIS.
Perhaps once these and other tilings
have been analyzed, it will be possible to state a 
more precise algorithmic version of the coloring book method.

We also thank the referees for their helpful comments.



\begin{table}
\caption{Calculation of coordination sequence for $3.12.12$ tiling. Top row gives  
distance $n$ from central vertex, rows indicate the eight types of vertex, last row is grand total.}\label{Tab31212}
$$
\begin{array}{c|cccccccc}
n: & 8k & 8k+1 & 8k+2 & 8k+3 & 8k+4 & 8k+5 & 8k+6 & 8k+7 \\
\hline
(i)   & 4 & 6 & 6 & 4 & 4 & 6 & 6 & 4 \\
(ii)  & 2 & 3 & 4 & 4 & 4 & 4 & 4 &  3 \\
(iii) & k & k & k & k & k & k & k+1 & k+1 \\
(iv) & 2k-1 & 2k-1 & 2k-1 & 2k-1 & 2k & 2k & 2k & 2k \\
(v) & 2k-1 & 2k-1 & k & k & 2k & 2k & k & k  \\
(vi) & k & k & k & k & k & k+1 & k+1 & k+1 \\
(vii) & 2k-1 & 2k-1 & 2k-1 & 2k & 2k & 2k & 2k & 2k+1 \\
(viii) & 2k-1 & k & k & 2k & 2k & k & k & 2k+1 \\
\hline
a(n) & 20k-2 & 18k+3 & 16k+6 & 18k+6 & 20k+8 & 18k+12 & 16k+14 & 18k+15  
\end{array}
$$
\end{table}

\begin{table}
\caption{Presentations for the groups of the eleven uniform tilings, giving
section number if mentioned above, number of generators $g$,
and corresponding sequence number. }\label{TabCox} 
$$
\begin{array}{lllll}
\mbox{Tiling} & $\S$ & g & \mbox{Presentation}  & \mbox{Sequence} \\
3^6       & &3& RST=1, RS=SR &  \seqnum{A008458} \\
3^4 \ldot6    & \ref{Sec33336} &3& R^2 = S^3 = T^6 = RST = 1 &  \seqnum{A250120} \\
3^3 \ldot 4^2   & &4& R^2=T^2=U^2=SUT=1, RS=SR &  \seqnum{A008706} \\
3^2 \ldot 4.3.4 & \ref{SecDual} &3& R^2= RST= (ST^{-1})^2 = 1 &  \seqnum{A219529} \\
3.4.6.4  & \ref{Sec3464}  &2& R^3=S^6=RSRS = 1 &   \seqnum{A008574} \\
3.6.3.6   & &2& R^3=S^3 = (RS)^3 = 1 &  \seqnum{A008579} \\
3.12^2    & \ref{Sec31212} &2& R^2=S^3 = (RS)^6 = 1 &  \seqnum{A250122} \\
4^4     & \ref{SecM}   &2& RS=SR &  \seqnum{A008574} \\
4.6.12   &  &3& R^2= S^2= T^2= (TS)^2= (RT)^3= (SR)^6 = 1 &  \seqnum{A072154} \\
4.8^2     & \ref{Sec488} &2& R^2=S^4=(RS)^4 = 1 &  \seqnum{A008576} \\
6^3      &  &3& R^2= S^2= T^2=1, RST=TSR &  \seqnum{A008486}\\
\end{array}
$$
\end{table}


\begin{references}

\reference{Baake, M. \& Grimm, U. (1997).
\emph{Z. Krist.}  \textbf{212}, 253--256.}

\reference{Bacher, R. \& de la Harpe, P. (2018).
\emph{Internat. Math. Res. Notices}, \textbf{2018:5}, 1532--1584.} 

\reference{Bacher, R., de la Harpe, P. \& Venkov, B. (1997). 
\emph{C. R. Acad. Sci. Paris, S\'{e}r. 1},  \textbf{325}, 1137--1142.}

\reference{Benson, M. (1983).
\emph{Invent. Math.}, \textbf{73}, 251--269.}

\reference{Blatov, V. A., Shevchenko, A. P. \& Proserpio, D. M. (2014).
\emph{Cryst. Growth Des.}, \textbf{14}, 3576Ð-3586.}

\reference{Bosma, W., Cannon, J. \& Playoust, C . (1997).
\emph{J. Symbolic Comput.}, \textbf{24}, 235Ð-265.}

\reference{Chavey, D. (1989).
\emph{Computers \& Mathematics with Applications}, \textbf{17:1--3}, 147--165.}

\reference{Conway, J. H., Burgiel, H. \& Goodman-Strauss, C. (2008).
\emph{The Symmetries of Things},
Wellesley, MA: A.~K.~Peters.}

\reference{Conway, J. H. \& Sloane, N. J. A. (1997).
\emph{Proc. Royal Soc. London, Series A},
\textbf{453}, 2369--2389.}

\reference{Coxeter, H. S. M. \& W. O. J. Moser (1984).
\emph{Generators and Relations for Discrete Groups}, 
4th. ed., Springer.}

\reference{de la Harpe, P. (2000). 
\emph{Topics in Geometric Group Theory}, Univ. Chicago Press,}

\reference{Eon, J.-G. (2002). 
\emph{Acta Cryst.}, \textbf{A58}, 47--53.}

\reference{Eon, J.-G. (2004). 
\emph{Acta Cryst.}, \textbf{A60}, 7--18.}

\reference{Eon, J.-G. (2007). 
\emph{Acta Cryst.}, \textbf{A63}, 53--65.}

\reference{Eon, J.-G. (2013). 
\emph{Acta Cryst.}, \textbf{A69}, 3 pages.}

\reference{Eon, J.-G. (2016). 
\emph{Acta Cryst.}, \textbf{A72}, 268--293.}

\reference{Eon, J.-G. (2018). 
\emph{Symmetry}, \textbf{10}, 13 pages.}

\reference{Epstein, D. B. A., Holt, D. F. \& Rees, S. E. (1991). 
\emph{J. Symbolic Computation}, \textbf{12}, 397--414.}

\reference{Galebach, B. (2018). 
\emph{N-uniform Tilings},
http://probabilitysports.com/tilings.html.}

\reference{Goodman-Strauss, C. (2009).
\emph{Theoret. Comp. Sci.}, \textbf {410}, 1534--1549.}

\reference{Grosse-Kunstleve, R. W., Brunner, G. O. \& Sloane, N. J. A. (1996).
\emph{Acta Cryst.}, \textbf{A52}, 879--889.}

\reference{Gr\"{u}nbaum, B. \& Shephard, G. C. (1977).
\emph{Math. Magazine}, \textbf{50}, 227--247.}

\reference{Gr\"{u}nbaum, B. \& Shephard, G. C. (1987).
\emph{Tilings and Patterns},
New York: W.~H.~Freeman.}

\reference{Johnson, D. L. (1976).
\emph{Presentations of Groups},
Cambridge.}

\reference{Knuth, D. E. \& Bendix. P. B. (1970). 
in J.~Leech, ed.,
\emph{Computational Problems in Abstract Algebra},
Oxford: Pergamon, pp. 263--297.}

\reference{O'Keeffe, M. (1995). 
\emph{Z. Krist.}, \textbf{210}, 905-908.}

\reference{O'Keeffe, M. \& Hyde, B. G. (1980). 
\emph{Phil. Trans. Royal Soc. London}, Series A, Mathematical and Physical Sciences, 
\textbf{295:1417}, 553Ð-618.}

\reference{O'Keeffe, M., Peskov, M. A., Ramsden, S. J. \& Yaghi, D. M. (2008).
\emph{Accounts of Chemical Research}, \textbf{41.12}, 1782--1789.}

\reference{The OEIS Foundation Inc. (2018).
\emph{The On-Line Encyclopedia of Integer Sequences},\\
https://oeis.org.}

\reference{Shutov, A. V. (2003).
\emph{Zap. Nauchn. Sem. S.-Peterburg. Otdel. Mat. Inst. Steklov. (POMI)},
\emph{Anal. Teor. Chisel i Teor. Funkts.}, \textbf{19}, 188--197, 203; 
English translation in \emph{J. Math. Sci. (N.Y.)},
\textbf{129}, 3922--3926.}

\reference{Zhuravlev, V. G. (2002).
\emph{St. Petersburg. Math. J.}, \textbf{13}, 201--220.}


\end{references}
\end{document}